\documentclass[11 pt]{amsart}

\usepackage{amsmath,amssymb,amsfonts,graphics}
\usepackage{verbatim}
\newcommand{\N}{\mathbb{N}}
\newcommand{\Z}{\mathbb{Z}}

\newcommand{\A}{\mathbb{A}}
\newcommand{\C}{\mathbb{C}}
\newcommand{\F}{\mathbb{F}}
\newcommand{\G}{\mathbb{G}}

\newcommand{\T}{\mathbb{T}}

\newcommand{\mc}{\mathcal}
\newcommand{\id}{\text{id}}

\newtheorem{thm}{Theorem}[section]
\newtheorem{cor}[thm]{Corollary}
\newtheorem{lem}[thm]{Lemma}
\newtheorem{prop}[thm]{Proposition}
\theoremstyle{definition}
\newtheorem{defn}[thm]{Definition}

\theoremstyle{remark}
\newtheorem{rem}[thm]{Remark}
\newtheorem{notat}[thm]{Notation}

\numberwithin{equation}{section}

\begin{document}

\title[Approximation Properties for Free Quantum Groups]
{Approximation properties for free orthogonal and free unitary quantum groups}

\date{\today}

\author{Michael Brannan}

\address{Michael Brannan: Department of Mathematics and Statistics, Queen's University, Kingston, Ontario, Canada, K7L 3N6.}
\email{mbrannan@mast.queensu.ca}

\keywords{Approximation properties, Haagerup approximation property, quantum groups, metric approximation property, free probability.}
\thanks{2010 \it{Mathematics Subject Classification:}
\rm{46L54, 46L10 (Primary); 46L65 (Secondary)}}

\begin{abstract}
In this paper, we study the structure of the reduced C$^\ast$-algebras and von Neumann algebras associated to the free orthogonal and free unitary quantum groups. We show that the reduced von Neumann algebras of these quantum groups always have the Haagerup approximation property.  Combining this result with a Haagerup-type inequality due to Vergnioux \cite{Ve}, we also show that the reduced C$^\ast$-algebras always have the metric approximation property.
\end{abstract}
\maketitle

\section{Introduction} \label{section_intro}

In \cite{Wa0, Wa}, S. Wang constructed a new family of compact quantum groups - the so-called free unitary and free orthogonal quantum groups - by considering certain highly noncommutative analogues of the C$^\ast$-bialgebras of continuous functions on the classical unitary and orthogonal groups $U_N$ and $O_N$.  More precisely, for any $N \ge 2$, the \textit{free unitary quantum group} (of dimension $N$) is the pair $U_N^+ := (A_u(N),\Delta)$, where $A_u(N)$ is the universal C$^\ast$-algebra defined by $N^2$ generators $\{u_{ij}:1 \le i,j \le N\}$, subject to the relations which make the matrices \begin{eqnarray} \label{eqn_def_A_U(N)} U := [u_{ij}]_{1 \le i,j \le N}  &\textrm{and}& \overline{U} :=  [u_{ij}^*]_{1 \le i,j \le N}
\end{eqnarray}   
unitary, and $\Delta:A_u(N) \to A_u(N) \otimes_{\textrm{min}} A_u(N)$ is the coproduct determined by 
\begin{eqnarray} \label{eqn_def_A_U(N)2} \Delta (u_{ij}) = \sum_{k=1}^N u_{ik}\otimes u_{kj} && (1 \le i,j \le N). 
\end{eqnarray}  
Similarly, the \textit{free orthogonal quantum group} is the pair $O_N^+ := (A_o(N), \Delta)$, where $A_o(N)$ is the universal C$^\ast$-algebra defined by $N^2$ self-adjoint generators $\{v_{ij}=v_{ij}^*:1 \le i,j \le N\}$, subject to the relations which make the matrix \begin{eqnarray} \label{eqn_defn_orth}
V := [v_{ij}]_{1 \le i,j \le N}
\end{eqnarray} unitary, and $\Delta:A_o(N) \to A_o(N) \otimes_{\textrm{min}} A_o(N)$ is the coproduct defined by 
\begin{eqnarray} \label{eqn_defn_orth2} \Delta (v_{ij}) = \sum_{k=1}^N v_{ik}\otimes v_{kj} && (1 \le i,j \le N). 
\end{eqnarray}  

Since the introduction of these quantum groups, the series $\{O_N^+\}_{N \ge 2}$ and $\{U_N^+\}_{N \ge 2}$ have been the object of intensive study from both quantum probabilistic and operator algebraic perspectives \cite{Ba0, Ba, BaCo, BaCuSp, BaCuSp2, CuSp, VaVe, Ve}.  On the probabilistic side, it is becoming increasingly apparent that the same role played by the groups $O_N$ and $U_N$ in classical probability theory, is played by $O_N^+$ and $U_N^+$ in Voiculescu's free probability theory \cite{BaCuSp, BaCuSp2, Cu}.  For example, Freedman's theorem characterizing conditionally independent Gaussian families in terms of $O_N$ and $U_N$ invariant distributions, has a free analogue when these groups are replaced by $O_N^+$ and $U_N^+$ \cite{Cu}.  The series $\{U_N^+\}_{N \ge 2}$ also provides canonical examples of ``quantum Haar unitary'' random matrices, which can be used to prove asymptotic \textit{operator valued} freeness results - a phenomenon which does not occur when using classical Haar unitary random matrices \cite{CuSp}.  From the perspective of operator algebras, the reduced C$^\ast$-algebras and von Neumann algebras arising from the GNS representations of the Haar traces on $O_N^+$ and $U_N^+$ form an interesting class of operator algebras, about which much is still not known.
  
For $\G = O_N^+$ or $U_N^+$, let $C(\G)$ denote the reduced C$^\ast$-algebra of $\G$, and let $L^\infty(\G) = C(\G)^{\prime \prime}$ denote the reduced von Neumann algebra of $\G$.  In \cite{Ba} and \cite{VaVe}, the algebras $C(\G)$ and $L^\infty(\G)$ were shown to share many structural properties with the reduced C$^\ast$- and von Neumann algebras of the free groups $\F_k$ ($k \ge 2$).  For example, $C(U_N^+)$ is always a non-nuclear, exact, and simple C$^\ast$-algebra, and $L^\infty(U_N^+)$ is always a solid, non-injective II$_1$-factor.  The same properties hold for $C(O_N^+)$ and $L^\infty(O_N^+)$ when $N \ge 3$. In \cite{Ba}, it is also shown that $L^\infty(U_2^+)$ is isomorphic to $L(\F_2)$, the von Neumann factor generated by the left regular representation of $\F_2$.  Given the above structural properties, it is natural to ask whether the algebras $L^\infty(\G)$ ($\G = O_N^+$ or $U_N^+,$ $N \ge 3$) also give rise to free group factors?  At this time, a direct answer to this question seems to be out of reach. Nevertheless, we can try and shed light on this question by asking what other structural properties of $C(\G)$ and $L^\infty(\G)$ are shared by the reduced free group algebras $C^*_\lambda(\F_k)$ and $L(\F_k)$ ($k \ge 2$)?    

In this paper, we consider this second question, and in particular look at what sorts of approximation properties these operator algebras have.  We prove that the von Neumann algebras $L^\infty(O_N^+)$ and $L^\infty(U_N^+)$ have the Haagerup approximation property for all $N \ge 2$, which answers questions posed by Vaes \cite{Va} and Vergnioux \cite[Introduction]{Ve1}.  Recall that a finite von Neumann algebra $(M, \tau)$ (with faithful normal trace $\tau:M \to \C$) has the \textit{Haagerup approximation property} (HAP) if there exists a net $\{\Phi_\lambda\}_{\lambda \in \Lambda}$ of normal, unital, completely positive, $\tau$-preserving maps on $M$ such that:
\begin{enumerate}
\item For each $\lambda \in \Lambda$, the unique $L^2$-extension $\hat{\Phi}_\lambda: L^2(M) \to L^2(M)$ is a compact operator, and  
\item For each $x \in M$, $\lim_{\lambda \in \Lambda}\|\hat{\Phi}_\lambda(\hat{x}) - \hat{x}\|_{L^2(M)} = 0,$ where $x \mapsto \hat{x}$ is the usual embedding $M \hookrightarrow L^2(M)$. 
\end{enumerate}  
Conditions $(1)$ and $(2)$ of course imply that the identity map $\id_{M \to M} = \lim_{\lambda \in \Lambda} \Phi_\lambda$ in the point $\sigma$-weak topology of $B(M)$.  Note that the HAP is a very important structural property for a finite von Neumann algebra, since it is an isomorphism invariant (see Remark \ref{rem_iso_inv}), and it is also a useful tool for proving other approximation properties, such as the weak$^\ast$-(completely) bounded approximation property (see \cite{RiHo, RiXu}).

In \cite{Ha}, Haagerup established the positive definiteness of the exponentiated length functions $\{g \mapsto e^{-|g|t}\}_{t \ge 0}$ on the free groups $\F_k$ ($k \ge 2$), which yields the HAP for the free group factors $L(\F_k)$.  Combining this result with the Haagerup inequality for the free groups \cite[Lemma 1.4]{Ha}, Haagerup also showed that the reduced C$^\ast$-algebras $C^*_\lambda(\F_k)$ have the metric approximation property (despite being non-nuclear).  Combining our proof of the HAP for $L^\infty(O_N^+)$ and $L^\infty(U_N^+)$ with a version of the Haagerup inequality for these quantum groups proved by Vergnioux \cite{Ve}, we are also able to show that the reduced C$^\ast$-algebras $C(O_N^+)$ and $C(U_N^+)$ have the metric approximation property.  This answers a question posed by Vergnioux in \cite[Introduction]{Ve}.   We also show that $C(U_2^+)$ has the (stronger) completely contractive approximation property (see Remark \ref{rem_CCAP_quest}).      

The remainder of this paper is organized as follows:  Section \ref{section_prelims} contains a brief review of the basic results on C$^\ast$-algebraic compact quantum groups that we will need.  In Section \ref{section_general_results}, we restate the definition of the Haagerup approximation property in the context of compact quantum groups of Kac type, and prove some general results about normal completely positive maps on their reduced von Neumann algebras.  In Section \ref{section_orthogonal}, we prove that $L^\infty(O_N^+)$ always has the HAP, and deduce the same result for $L^\infty(U_N^+)$ using a free product representation of $L^\infty(U_N^+)$ due to Banica \cite{Ba}.  The two main ingredients to our proof of the HAP for $L^\infty(O_N^+)$ are $(1)$ an averaging result for state-induced normal completely positive maps on compact Kac algebras (Theorem \ref{thm_central_cp_maps}), and $(2)$ the fact that the spectral measure of the fundamental character of $O_N^+$ (relative to the Haar trace) is always Wigner's semicircle law (Theorem \ref{thm_banica_orthogonal}).  In Section \ref{section_MAP}, we use Vergnioux's Haagerup inequalities for $O_N^+$ and $U_N^+$, and the results of Section \ref{section_orthogonal}, to prove that the C$^\ast$-algebras $C(O_N^+)$ and $C(U_N^+)$ have the metric approximation property.  We close with Remark \ref{rem_bai_L1}, where we use the results of Section \ref{section_MAP} to show that the Banach algebra predual $L^1(\G):=L^\infty(\G)_*$ ($\G = U_N^+$ or $O_N^+$) always has a central approximate identity which is bounded in the multiplier norm on $L^1(\G)$.  

\subsection*{Acknowledgements}  The author thanks his doctoral supervisors James A. Mingo and Roland Speicher for useful discussions.  The writing of this paper was completed while the author participated in the ``Bialgebras in Free Probability'' programme at the Erwin Schr\"odinger International Institute for Mathematical Physics.  The author thanks the Institute for their kind hospitality and the excellent working environment.  This research was partially supported by an NSERC Canada Graduate Scholarship.  

\section{Preliminaries and Notation} \label{section_prelims}

In this section we review the basic theory of compact quantum groups, as developed by Woronowicz in \cite{Wo2} (see also the excellent book \cite{Ti}).  For the remainder of this paper, we will use the symbol $\otimes$ when writing the minimal tensor product of a pair of C$^\ast$-algebras.  Similarly, $\overline{\otimes}$ will be used to denote the von Neumann tensor product, and $\otimes_{\textrm{alg}}$ will denote the purely algebraic tensor product of two algebras.  A basic familiarity with the theory of completely bounded maps will be assumed.  In particular, $B(X)$ (resp. $CB(X)$) will always denote the algebra of bounded (resp. completely bounded) operators on a Banach (resp. operator) space $X$.  Our reference for this will be \cite{Pa, Pi}. 

A \textit{compact quantum group} (CQG) is a pair $\G = (A, \Delta)$ where
\begin{itemize}
\item $A$ is a unital C$^\ast$-algebra,
\item $\Delta:A \to A \otimes A$ is a unital $\ast$-homomorphism satisfying the coassociativity relation 
\begin{equation}
(\id_A \otimes \Delta) \circ \Delta = (\Delta \otimes \id_A) \circ \Delta, 
\end{equation}
and 
\item $\G$ satisfies the \textit{cancellation property}.  That is, the sets \begin{eqnarray*} \Delta(A)(1_A \otimes A), &\textrm{and}& \Delta(A)(A \otimes 1_A),
\end{eqnarray*} are linearly dense in $A\otimes A$. 
\end{itemize}

From these three axioms, it follows that any CQG $\G = (A, \Delta)$ admits a unique state $h:A \to \C$, called the \textit{Haar state}, which satisfies the bi-invariance condition 
\begin{eqnarray}
(h \otimes \id_A)\Delta(x) = (\id_A \otimes h) \Delta(x) = h(x)1_A && (x \in A).
\end{eqnarray}
We denote by $L^2(\G)$ the usual GNS Hilbert space obtained from equipping $A$ with the sesquilinear form $(x,y) \mapsto h(y^*x)$, and denote by $\pi_h:A \to B(L^2(\G))$ the GNS representation.  We also put $C(\G) = \pi_h(A) \subseteq B(L^2(\G))$, and $L^\infty(\G) = C(\G)^{\prime \prime} \subseteq B(L^2(\G))$.  $C(\G)$ is called the \textit{reduced C$^\ast$-algebra} of $\G$, and $L^\infty(\G)$ is called the \textit{reduced von Neumann algebra} of $\G$.  Since $h = (h \otimes h) \circ \Delta$, it follows that the coproduct $\Delta$ determines a coproduct $\Delta_r:C(\G) \to C(\G) \otimes C(\G)$ satisfying $\Delta_r \circ \pi_h = (\pi_h \otimes \pi_h) \circ \Delta$.  Furthermore, $\Delta_r$ extends to a normal $\ast$-homomorphism $\Delta_r:L^\infty(\G) \to L^\infty(\G) \overline{\otimes} L^\infty(\G).$  The CQG $\G_r = (C(\G), \Delta_r)$ is called the reduced version of $\G$, and $(L^\infty(\G), \Delta_r)$ is a von Neumann algebraic CQG, in the sense of \cite{KuVa}.  The Haar state $h_r:C(\G) \to \C$ is given by $h_r (\pi_h(x)) = h(x)$ for each $x \in A$, and extends to a faithful normal state on $L^\infty(\G)$.  From now on we will drop the subscript $r$ from $h_r$, and just write $h$ for the Haar state in all possible situations.  No confusion should arise from this.

For any CQG $\G = (A, \Delta)$, the Banach space $A^*$ of continuous linear functionals on $A$ comes equipped with the associative convolution product $$(\varphi, \psi) \mapsto \varphi * \psi =  (\varphi \otimes \psi) \circ \Delta.$$  With this product, $A^*$ is a (completely contractive) Banach algebra.  The predual $L^1(\G) := L^\infty(\G)_*$ of the reduced von Neumann algebra of $\G$ can be completely isometrically identified with a closed subspace of $A^*$ via the dual pairing 
\begin{eqnarray} \label{eqn_dual_pairing}
\langle x, \omega \rangle = \langle \omega, \pi_h(x) \rangle && (x \in A, \omega \in L^1(\G)). 
\end{eqnarray} 
With this identification, $L^1(\G)$ is a closed two-sided ideal in $A^*$ (see \cite{Da}, \cite[Page 913-914]{KuVa1}).

\begin{defn} \label{defn_corep}
An ($n$-dimensional) \textit{corepresentation} of a CQG $\G = (A, \Delta)$ is a matrix $U = [u_{ij}]_{1 \le i,j \le n} \in M_n(A) = M_n(\C) \otimes A$ with the property that \begin{eqnarray*}
\Delta(u_{ij}) = \sum_{k=1}^n u_{ik} \otimes u_{kj}  && (1 \le i,j \le n).
\end{eqnarray*}  
We say that $U$ is a \textit{unitary corepresentation} if in addition $U$ is a unitary element of $M_n(A)$.
\end{defn}

If $U = [u_{ij}]_{1 \le i,j \le n}$, is a corepresentation, the matrix $\overline{U} = [u_{ij}^*]_{1 \le i,j \le n}$ is also a corepresentation, called the \textit{conjugate} of $U$.  We note that $\overline{U}$ may not be a unitary corepresentation, even if $U$ is unitary. 

If $U \in M_n(\C) \otimes A$ and $V \in M_m(\C) \otimes A$ are two corepresentations of $\G$, we define the vector space $$\textrm{Hom}(U,V) = \{T \in B(\C^n, \C^m): \ (T \otimes 1_A)U = V(T \otimes 1_A)\}.$$  An element $T \in \textrm{Hom}(U,V)$ is called an \textit{intertwiner} from $U$ to $V$.  The corepresentation $U$ is called \textit{irreducible} if $\textrm{Hom} (U,U) = \C \id$, and $U$ is irreducible if and only if $\overline{U}$ is.  $U$ and $V$ are called \textit{(unitarily) equivalent} corepresentations if there exists an invertible (unitary) operator $T \in \textrm{Hom}(U,V)$, and we write $U \cong V$. 

The matrix $$U \boxtimes V := [u_{ij}v_{kl}]_{\substack{1 \le i,j \le n \\ 1 \le k,l \le m}} \in M_n(\C) \otimes M_m(\C) \otimes A,$$ is again a corepresentation of $\G$, called the \textit{tensor product} of $U$ and $V$.

One can of course define the notion of an infinite dimensional (unitary) corepresentation of $\G$ (see \cite[Section 2]{Wo2}), however we will not need this generality here, mainly due to the following theorem:

\begin{thm} \label{thm_reducibility} (\cite{Wo2}) 
Every irreducible corepresentation of a CQG is finite dimensional and equivalent to a unitary one. Furthermore, every
unitary corepresentation is unitarily equivalent to a direct sum of irreducibles.
\end{thm}

Let $\{U^\alpha = [u^\alpha_{ij}]_{1 \le i,j \le d_\alpha} : \alpha \in \A \}$ be a maximal family of finite-dimensional irreducible unitary
corepresentations of $\G = (A, \Delta)$, with $\alpha_0$ being the index corresponding to the trivial corepresentation $1_A \in A$, and $U^{\overline{\alpha}}$ denoting the representative of the class of $\overline{U^\alpha}$.  Let $\mc A$ denote the subspace of $A$ spanned by the matrix elements $\{u^\alpha_{ij}: \ 1 \le i,j \le d_\alpha, \alpha \in \A \}$. Then $\mc A$ is a Hopf $\ast$-algebra with faithful Haar state given by $h|_\mc A$, $\mathcal A$ is norm dense in $A$, and the set $\{u^\alpha_{ij}: \ 1 \le i,j \le d_\alpha, \alpha \in \A \}$ is a linear basis for $\mc A$.  The coproduct $\Delta_\mc A:\mc A \to \mc A \otimes_{\textrm{alg}} \mc A $ is just the restriction $\Delta|_{\mc A}$,  the \textit{coinverse} $\kappa: \mc A \to \mc A$ is the antihomomorphism given by 
\begin{eqnarray*}
\kappa(u_{ij}^\alpha)= (u_{ji}^\alpha)^*  && (1 \le i,j \le d_\alpha, \ \alpha \in \A),
\end{eqnarray*}
and the \textit{counit} $\epsilon:\mc A \to \C$ is the $\ast$-character given by 
\begin{eqnarray*}
\epsilon(u_{ij}^\alpha) = \delta_{ij} && (1 \le i,j \le d_\alpha, \ \alpha \in \A).
\end{eqnarray*}
The Hopf $\ast$-algebra $\mc A$ is unique in the sense that if $\mathcal B \subseteq A$ is any other dense Hopf $\ast$-subalgebra, then $\mathcal B = \mathcal A$ (see \cite[Theorem 5.1]{BeMuTu0}).  

In general, the counit $\epsilon:\mc A \to \C$ and the coinverse $\kappa: \mc A \to \mc A$ defined above cannot be extended to bounded linear maps on $A$.  A CQG $\G = (A, \Delta)$ is called \textit{co-amenable} if the map $\pi_h(x) \mapsto \epsilon(x)$ ($x \in \mc A$) extends to a character of $C(\G)$.   In this paper we will mainly deal with CQGs of Kac type, which are precisely those for which the coinverse $\kappa$ can be extended to a bounded map on $A$.  

\begin{defn} \label{defn_Kac}
A CQG $\G = (A, \Delta)$ is said to be of \textit{Kac type} if any one of the following equivalent conditions is satisfied.
\begin{enumerate}
\item $\kappa: \mathcal A \to \mc A$ has a continuous extension to a $\ast$-antihomomorphism $\kappa:A \to A$.
\item $\kappa^2 = \id_\mc A$. 
\item The Haar state $h:A \to \C$ is a trace.
\end{enumerate}
\end{defn}  
See \cite{BaSk} for proofs of the above equivalences.  From the universality of the C$^\ast$-algebras of $A_o(N)$ and $A_u(N)$ defined in Section \ref{section_intro}, it is not difficult to see that condition $(1)$ of Definition \ref{defn_Kac} is satisfied for these algebras.  So $O_N^+$ and $U_N^+$ are of Kac type.  Also, $O_N^+$ is not co-amenable for all $N \ge 3$, and $U_N^+$ is not co-amenable for all $N \ge 2$ \cite{Ba}. 

When $\G = (A, \Delta)$ is of Kac type and $U^\alpha$ ($\alpha \in \A$) is an irreducible unitary corepresentation (using the above notation), it is readily checked that $\overline{U^\alpha}$ is also unitary.  We can therefore assume from now on that our representatives $\{U^\alpha: \alpha \in \A\}$ have been chosen so that $U^{\overline{\alpha}} = \overline{U^\alpha}$.  In this case, we have
\begin{eqnarray} \label{eqn_kappa_action}
\kappa(u^\alpha_{ij}) = (u^\alpha_{ji})^* = u^{\overline{\alpha}}_{ji} && (1 \le i,j \le d_\alpha, \ \alpha \in \A).
\end{eqnarray}
Furthermore, the family of matrix elements $\{\sqrt{d_\alpha}u_{ij}^\alpha: \alpha \in \A, 1 \le i,j \le d_\alpha\} \subset \mc A$ always forms an orthonormal basis for $L^2(\G)$ in the Kac algebraic setting.

Finally, when $\G = (A,\Delta)$ is of Kac type, the map $\kappa_r:C(\G) \to C(\G)$ defined by $\kappa_r(\pi_h(x)) = \pi_h(\kappa(x))$ $(x \in A)$, defines the coinverse for the reduced quantum group $\G_r$.   Furthermore, $\kappa_r$ extends to a normal $\ast$-antiautomorphism $\kappa_r:L^\infty(\G) \to L^\infty(\G)$.

\section{The Haagerup Approximation Property for Compact Quantum Groups of Kac Type} \label{section_general_results}

In this section, we discuss the Haagerup approximation property in the context of CQGs $\G = (A,\Delta)$ of Kac type, and prove some general results on constructing normal completely positive maps on $L^\infty(\G)$ from states on $A$.   

We will say that a CQG $\G = (A, \Delta)$ of Kac type has the \textit{Haagerup approximation property (HAP)} if the finite von Neumann algebra $(L^\infty(\G), h)$ has the HAP, as defined in Section \ref{section_intro}.  

\begin{rem} \label{rem_iso_inv} At this point, it is worthwhile to mention that the HAP for a finite von Neumann algebra $(M, \tau)$ does not actually depend on the particular choice of trace $\tau$.  I.e., if $\tau^\prime:M \to \C$ is another faithful, normal, finite trace, then $(M,\tau^\prime)$ also has the HAP.  Furthermore, in the definition of the HAP, the requirement that the maps $\Phi_\lambda:M \to M$ be $\tau$-preserving can be relaxed to $\tau \circ \Phi_\lambda \le \tau$, and this yields the same class of von Neumann algebras \cite{Jo1}.  From the trace-independence of the HAP, it is easily deduced that the HAP is an invariant for the isomorphism class of a finite von Neumann algebra.
\end{rem}

\begin{notat} Let $A$ be a C$^\ast$-algebra and let $T:A \to A$ be a linear map. We will say that $T$ is UCP when $T$ is unital and completely positive.  If $A$ is a von Neumann algebra, we will say that $T$ is NUCP if $T$ is normal, unital, and completely positive.
\end{notat}

In order to establish that a given CQG $\G = (A, \Delta)$ of Kac type has the HAP, we need to construct a net of NUCP, $h$-preserving maps on $L^\infty(\G)$ which satisfies certain properties.  In the classical case (when $\G$ corresponds either to a classical compact group $G$ or the dual of a discrete group $\Gamma$), a convenient place to look for such maps is within the set of convolution operators induced by states on the C$^\ast$-algebra $A$.  In the general situation, the same is true: states on $A$ yield concrete examples of NUCP maps on $L^\infty(\G)$.  We will now briefly outline this well known procedure.  
\begin{defn} \label{defn_conv}
Let $\G = (A, \Delta)$ be a CQG and $\varphi \in A^*$ a state.  The UCP $h$-preserving map $C_\varphi:A \to A$  given by  \begin{eqnarray*} C_\varphi x  = (\varphi \otimes \id_A) \Delta(x) && (x \in A),
\end{eqnarray*} 
is called the \textit{(left) convolution operator} associated to $\varphi$.
\end{defn}

Let $\pi_h:A \to C(\G)$ denote the GNS representation of the Haar state $h:A \to \C$.  The following lemma shows that the convolution operator $C_\varphi$ factors through the quotient $C(\G) \cong A/\ker \pi_h,$ and furthermore extends to a NUCP map on $L^\infty(\G)$.  This result is already known, and can be nicely stated (in a more general context) in terms of completely bounded multipliers of locally compact quantum groups (see \cite[Proposition 8.3]{Da} and \cite[Section 6]{Da1}).  In order to avoid dealing with Pontryagin duality and to keep this paper self-contained, we provide a proof of the following lemma. 

\begin{lem} \label{lem_CP_maps_states}
Let $\G = (A, \Delta)$ be a CQG, and let $\varphi \in A^*$ be a state.  Then there exists a unique UCP $h$-preserving map $S_\varphi:C(\G) \to C(\G)$ defined by \begin{eqnarray*}
S_\varphi(\pi_h(x)) = \pi_h(C_\varphi x) && (x \in A),  
\end{eqnarray*}
where $C_\varphi$ is the convolution operator given in Definition \ref{defn_conv}.  Furthermore, $S_\varphi$ extends uniquely to a NUCP $h$-preserving map 
$$S_\varphi:L^\infty(\G) \to L^\infty(\G).$$
\end{lem}

\begin{proof}
We claim that $C_\varphi ( \ker \pi_h) \subseteq \ker \pi_h$.  To see this, fix $x \in \ker \pi_h$.  Since $h$ is a KMS-state \cite[Example 8.1.22]{Ti}, it follows that $\ker \pi_h = \{x \in A : h(x^*x) = 0\}$.  By the Cauchy-Schwarz inequality for completely positive maps (\cite[Proposition 3.3]{Pa}), we have $$(C_\varphi x)^*(C_\varphi x) \le C_\varphi(x^*x).$$  Applying $h$ to this inequality and using the fact that $C_\varphi$ is $h$-preserving, we get $$h((C_\varphi x)^*(C_\varphi x) )  \le h(C_\varphi (x^*x)) = h(x^*x) = 0.$$  So $C_\varphi x \in \ker \pi_h$, proving the claim.  

Since $C_\varphi(\ker \pi_h) \subseteq \ker \pi_h,$ we may define a linear map $$S_\varphi:C(\G) \to C(\G),$$ by setting $$S_\varphi(\pi_h(x)) = \pi_h(C_\varphi x).$$  The fact that $S_\varphi$ is UCP and $h$-preserving follows from this formula and the fact that $C_\varphi$ has these properties.  

We now prove that $S_\varphi$ extends to a normal map on $L^\infty(\G)$.  Since $L^1(\G) \subseteq A^*$ is an ideal, we can consider the map $S_{\varphi \ast} \in CB(L^1(\G))$ given by \begin{eqnarray*} S_{\varphi \ast} (\omega) = \varphi * \omega && (\omega \in L^1(\G)).
\end{eqnarray*}
Then the adjoint $(S_{\varphi *})^* \in CB(L^\infty(\G))$ is normal, and $(S_{\varphi *})^*|_{C(\G)} = S_\varphi$.  Indeed, taking $x \in A$, $\omega \in L^1(\G)$, and using the identification (\ref{eqn_dual_pairing}), we have \begin{eqnarray*}
&& \langle \omega, S_\varphi \pi_h(x) \rangle = \langle \omega, \pi_h((\varphi \otimes \id_A) \Delta(x)) \rangle = \langle (\varphi \otimes \id_A) \Delta(x), \omega \rangle \\
&=& (\varphi \otimes \omega) \Delta(x)  = \langle x,  \varphi *\omega\rangle = \langle S_{\varphi *}\omega, \pi_h(x) \rangle \\
&=& \langle \omega, (S_{\varphi \ast})^* \pi_h(x) \rangle.
\end{eqnarray*}
Therefore $(S_{\varphi *})^*|_{C(\G)} = S_\varphi$.  Since $C(\G)$ is $\sigma$-weakly dense in $L^\infty(\G)$, we find that $(S_{\varphi *})^*$ is the unique normal extension of $S_\varphi$.  The fact that it is UCP and $h$-preserving follows from these properties for $S_\varphi$. 
\end{proof}

\begin{rem} \label{rem_action_S}
In order to keep our notation simple, we will identify the Hopf $\ast$-algebra $\mc A \subseteq A$ associated to a CQG $\G = (A,\Delta)$ with its image $\pi_h(\mathcal A) \subset C(\G)$ under the GNS representation $\pi_h:A \to C(\G)$.  This is possible because the Haar state $h|_{\mc A}$ is always faithful (and therefore $\pi_h|_{\mc A}$ is injective).  With this identification, the NUCP map $S_\varphi \in CB(L^\infty(\G))$ considered in Lemma \ref{lem_CP_maps_states} is given on the basis $\{u_{ij}^\alpha: 1 \le i,j \le d_\alpha, \  \alpha \in \A\} \subset \mc A$ by the formula 
\begin{eqnarray} \label{eqn_action_S}
S_\varphi(u_{ij}^\alpha) = (\varphi \otimes \id_\mathcal A)\Delta(u_{ij}^\alpha) = \sum_{k=1}^{d_\alpha} \varphi(u_{ik}^\alpha)u_{kj}^\alpha.  
\end{eqnarray}
We will also identify $C(\G)$ and $L^\infty(\G)$ with their dense images in $L^2(\G)$, via the GNS map $x \mapsto \hat{x}$.  Consequently, we have the natural inclusions $\mathcal A \subseteq C(\G) \subseteq L^\infty(\G) \subseteq L^2(\G),$ and $A = \overline{\mc A}^{\|\cdot\|_A}$.
\end{rem}

Equation (\ref{eqn_action_S}) gives a useful description of the NUCP map $S_\varphi:L^\infty(\G) \to L^\infty(\G)$ in terms of the action of the state $\varphi$ on $\mathcal A$.  Unfortunately, for a general CQG $\G = (A, \Delta)$, it is a difficult task to explicitly construct (nontrivial) states $\varphi \in A^*$.  On the positive side, when $\G$ is of Kac type, $S_\varphi$ can be ``averaged'' using the Haar trace to get a new NUCP map $T_\varphi$, which now only depends on the restriction of $\varphi$ to the subalgebra of $A$ generated by the irreducible characters of $\G$.   In order to introduce this averaging result, we first need some notation.   

\begin{notat} \label{notat_L2_alpha}
For each $\alpha \in \A$, denote by $L^2_\alpha(\G) \subset L^2(\G)$ the subspace spanned by the matrix elements $\{u_{ij}^\alpha: 1 \le i,j \le d_\alpha\}$ of the corepresentation $U^\alpha$.  Then $L^2(\G) = \ell^2 -  \bigoplus_{\alpha \in \A} L^2_\alpha(\G)$.  Denote by $p_\alpha:L^2(\G) \to L^2_\alpha(\G)$ the orthogonal projection, and let $$\chi_\alpha = (Tr_{d_\alpha} \otimes \id_\mc A)(U^\alpha) = \sum_{i=1}^{d_\alpha} u_{ii}^\alpha,$$ denote the \textit{irreducible character} of the corepresentation $U^\alpha$.  
\end{notat}

The following theorem is the main result of this section, and will be used to prove that $O_N^+$ has the HAP in Section \ref{section_orthogonal}. 

\begin{thm} \label{thm_central_cp_maps}
Let $\G = (A, \Delta)$ be a CQG of Kac type, and consider the unital C$^\ast$-subalgebra $\mathcal B = C^*\langle \chi_\alpha: \alpha \in \A \rangle \subseteq A$ generated by the irreducible characters $\{\chi_\alpha: \alpha \in \A\}.$  Then for any state $\psi \in \mc B^*$, 
\begin{enumerate}
\item The map \begin{eqnarray*}\hat{T}_\psi = \sum_{\alpha \in \A} \frac{\psi(\chi_{\overline{\alpha}})}{d_\alpha} p_\alpha, 
\end{eqnarray*}
is a unital contraction on $L^2(\G)$, and
\item The restriction $T_\psi = \hat{T}_\psi|_{L^\infty(\G)}$ defines a NUCP $h$-preserving map $T_\psi \in CB(L^\infty(\G))$.  Furthermore $T_\psi(C(\G)) \subseteq C(\G)$.
\end{enumerate}
\end{thm}

\begin{proof} (Of Theorem \ref{thm_central_cp_maps}(1).)
Since each $U^\alpha \in M_{d_\alpha}(A)$ is a unitary operator, $\|u_{ij}^\alpha\|_{A} \le 1$ for each $1 \le i,j \le d_\alpha$.  Therefore $$\Big| \frac{\psi(\chi_\alpha)}{d_\alpha} \Big|  \le d_\alpha^{-1}\Big\|\sum_{i=1}^{d_\alpha}u^\alpha_{ii}\Big\|_A \le d_{\alpha}^{-1}\sum_{i=1}^{d_\alpha}\|u^\alpha_{ii}\|_A \le 1.$$
Since the family of projections $\{p_\alpha\}_{\alpha \in \A} \subset B(L^2(\G))$ is orthogonal, $\hat{T}_\psi = \sum_{\alpha \in \A} \frac{\psi(\chi_{\overline{\alpha}})}{d_\alpha} p_\alpha$ satisfies $\|\hat{T}_\psi\|_{B(L^2(\G))} = \sup_{\alpha \in \A} \frac{|\psi(\chi_{\overline{\alpha}})|}{d_\alpha} \le 1$.  Since $\psi$ is a state, $\hat{T}_\psi 1_{L^\infty(\G)} = \psi(1_A)1_{L^\infty(\G)} = 1_{L^\infty(\G)}$. 
\end{proof}

The proof of Theorem \ref{thm_central_cp_maps}(2) is much more in involved.  The basic idea is to obtain $T_\psi$ by averaging (with respect to the Haar trace) the NUCP map $S_\varphi\in CB(L^\infty(\G))$, where $\varphi \in A^*$ is any Hahn-Banach extension of $\psi$. We will first need a few preparatory results.  

\begin{lem} \label{lem_conj_by_kappa}
Let $\G$ be a CQG of Kac type and let $S \in CB(L^\infty(\G))$ be unital, normal and $h$-preserving.  Then the map $$\kappa_r \circ S  \circ \kappa_r:L^\infty(\G) \to \L^\infty(\G)$$ is also unital, normal and $h$-preserving, with $\|\kappa_r \circ S  \circ \kappa_r\|_{cb} = \|S\|_{cb}.$  Furthermore,  we have $S(C(\G)) \subseteq C(\G)$ if and only if $(\kappa_r \circ S  \circ \kappa_r)C(\G) \subseteq C(\G)$.  
\end{lem}

\begin{proof}
This is actually a general result for Kac algebras, proved in \cite[Proposition 2.4]{KrRu} for example.  We include a proof for the reader's convenience.  

First of all,  $S$ and $\kappa_r$ are unital, normal, $h$-preserving, and $\kappa_r(C(\G))  = C(\G)$, it follows that $\kappa_r \circ S  \circ \kappa_r$ is also unital, normal and $h$-preserving.  Furthermore, $S(C(\G)) \subseteq C(\G)$ if and only if $(\kappa_r \circ S  \circ \kappa_r)C(\G) \subseteq C(\G)$.   Since $\kappa_r$ is a $\ast$-antiautomorphism, the map $$L^\infty(\G) \owns x \mapsto R(x) =  \kappa_r(x)^*,$$ defines a conjugate-linear $\ast$-automorphism of $L^\infty(\G)$, and in particular is a complete isometry.  This means that for any matrix $[x_{ij}] \in M_n(L^\infty(\G))$, \begin{eqnarray*}
&& \|[x_{ij}]\| = \|[x_{ij}]^*\| = \|[x_{ji}^*]\| = \|[\kappa_r(\kappa_r(x_{ji}^*))]\| = \|[\kappa_r(\kappa_r(x_{ji}))^*]\| \\
&=& \|(\id_{M_n(\C)} \otimes R)[\kappa_r (x_{ji})]\| = \|[\kappa_r (x_{ji})]\|.
\end{eqnarray*}
Using the above equality of matrix norms, we have for any matrix $[y_{ij}] \in M_n(L^\infty(\G))$, 
\begin{eqnarray*}
&& \|(\kappa_r \circ S \circ \kappa_r)^{(n)}[y_{ij}]\| = \|[\kappa_r (S \kappa_r(y_{ij}))]\| = \|[S(\kappa_r(y_{ji}))]\| \\
&\le& \|S\|_{cb} \|[\kappa_r(y_{ji})]\|  = \|S\|_{cb} \|[y_{ij}]\|. 
\end{eqnarray*}
Since $n \in \N$ was arbitrary, $\|\kappa_r \circ S \circ \kappa_r\|_{cb} \le \|S\|_{cb}$.  The reverse inequality follows from this one, since we may write $S = \kappa_r \circ (\kappa_r \circ S \circ \kappa_r) \circ \kappa_r.$
\end{proof}

Let $\G$ be a CQG of Kac type.  Since $(L^\infty(\G), h)$ is a finite von Neumann algebra and $$\Delta_r:L^\infty(G) \to L^\infty(\G) \overline{\otimes} L^\infty(\G)$$ is an injective normal $\ast$-homomorphism, it follows that $\Delta_r(L^\infty(\G))$ is a von Neumann subalgebra of $L^\infty(\G) \overline{\otimes} L^\infty(\G)$.  Therefore there exists a unique, normal, $h \otimes h$-preserving conditional expectation $$E:L^\infty(\G) \overline{\otimes}L^\infty(\G) \to \Delta_r(L^\infty(\G)).$$   The following lemma describes how this conditional expectation acts on $L^\infty(\G)\overline{\otimes}L^\infty(\G)$.

\begin{lem} \label{lem_action_cond_exp}
The conditional expectation $E:L^\infty(\G) \overline{\otimes}L^\infty(\G) \to \Delta_r(L^\infty(\G))$ is determined by the equations 
\begin{eqnarray*}
E(u_{ij}^\alpha \otimes u_{kl}^\beta) = \frac{\delta_{\alpha,\beta}\delta_{j,k}}{d_\alpha} \Delta_r(u_{il}^\alpha) && (\alpha, \beta \in \A, \ 1 \le i,j \le d_\alpha, \ 1 \le k,l \le d_\beta).
\end{eqnarray*}
Furthermore, $E(C(\G) \otimes C(\G)) \subset \Delta_r(C(\G))$.
\end{lem}

\begin{proof} As usual, we identify $L^\infty(\G)$ as a subspace of $L^2(\G)$ (and likewise we identify $L^\infty(\G) \overline{\otimes}L^\infty(\G)$ as a subspace of $L^2(\G)\otimes^2 L^2(\G)$).  With these identifications, $E$ is is realized simply as the orthogonal projection from $L^2(\G) \otimes^2 L^2(\G)$ onto $L^2(\Delta_r(L^\infty(\G))$. 

Recall that $\Delta_r:L^\infty(\G) \to \Delta_r(L^\infty(\G))$ is trace preserving, and therefore extends to an isometry with dense range for the $L^2$-norms (i.e. a unitary operator from $L^2(\G)$ to $L^2(\Delta_r(L^\infty(\G))$.)  Therefore $\Delta_r$ maps the orthonormal basis $\{d_{\alpha}^{1/2}u_{ij}^\alpha: 1 \le i,j \le d_\alpha, \alpha \in \A\}$ for $L^2(\G)$ to the orthonormal basis $\{\Delta_r(d_{\alpha}^{1/2}u_{ij}^\alpha): 1 \le i,j \le d_\alpha, \alpha \in \A\}$ for $L^2(\Delta_r(L^\infty(\G))$. Therefore we can write $E(u_{ij}^\alpha \otimes u_{kl}^\beta)$ as the $L^2$-convergent series
\begin{eqnarray*}
E(u_{ij}^\alpha \otimes u_{kl}^\beta) &=& \sum_{\gamma \in \A} \sum_{1 \le t,s \le d_{\gamma}} \langle u_{ij}^\alpha \otimes u_{kl}^\beta | \Delta_r(d_{\gamma}^{1/2}u_{ts}^\gamma)  \rangle \Delta_r(d_{\gamma}^{1/2}u_{ts}^\gamma) \\
&=& \sum_{\gamma \in \A} d_\gamma \sum_{1 \le t,s \le d_{\gamma}} \sum_{m=1}^{d_\gamma} \langle u_{ij}^\alpha \otimes u_{kl}^\beta | u_{tm}^\gamma \otimes u_{ms}^\gamma \rangle \Delta_r(u_{ts}^\gamma) \\
&=& \sum_{\gamma \in \A} d_\gamma \sum_{1 \le t,s \le d_{\gamma}} \sum_{m=1}^{d_\gamma} \frac{\delta_{i,t}\delta_{j,m}\delta_{k,m}\delta_{l,s}\delta_{\alpha, \gamma} \delta_{\beta, \gamma}}{d_\gamma^2} \Delta_r(u_{ts}^\gamma) \\
&=& \frac{\delta_{\alpha,\beta}\delta_{j,k}}{d_\alpha} \Delta_r(u_{il}^\alpha).
\end{eqnarray*}  
To see that $E(C(\G) \otimes C(\G)) \subseteq \Delta_r(C(\G))$, note that from the above formula for $E$, we have $E(\mc A \otimes_{\textrm{alg}}\mc A) \subseteq \Delta_r(\mc A)$.  Since $\mc A \otimes_{\textrm{alg}}\mc A$ is dense in $C(\G) \otimes C(\G)$, $\Delta_r(\mc A)$ is dense in $\Delta_r(C(\G))$ and $E$ is continuous, the result follows. 
\end{proof} 
We are now ready to prove Theorem \ref{thm_central_cp_maps}(2).

\subsection{Proof of Theorem \ref{thm_central_cp_maps}(2).}
Let $\varphi \in A^*$ be any Hahn-Banach extension of the given state $\psi \in \mc B^*$.  Since $\varphi(1_A) = \psi(1_A) = 1$, and $\|\varphi\| = \|\psi\| = 1$, $\varphi$ is also a state.  Let $S_\varphi \in CB(L^\infty(\G))$ be the NUCP $h$-preserving map defined in Lemma \ref{lem_CP_maps_states}, and consider the composition 
\begin{eqnarray} \label{eqn_averaging}
Q_\varphi = \Delta_r^{-1} \circ E \circ((\kappa_r \circ S_\varphi \circ \kappa_r) \otimes \id_{L^\infty(\G)}) \circ \Delta_r.
\end{eqnarray}
From Lemma \ref{lem_conj_by_kappa} above and the various definitions of the structure maps, it directly follows that $(\kappa_r \circ S_\varphi \circ \kappa_r) \otimes \id_{L^\infty(\G)}$, $\Delta_r$, $E$, and $\Delta_r^{-1}$ are all NUCP and trace-preserving maps.  Therefore their composition $Q_\varphi$ is also NUCP and $h$-preserving, and $Q_\varphi(C(\G)) \subseteq C(\G)$ by Lemmas \ref{lem_CP_maps_states}, \ref{lem_conj_by_kappa} and  \ref{lem_action_cond_exp}.  We now claim that $T_\psi = Q_\varphi$.  This will complete the proof.

To show this, it suffices by the normality of the composition $Q_\varphi$ to check equality on the basis $\{u_{ij}^\alpha: \alpha \in \A, \ 1 \le i,j \le d_\alpha\}$ of the the $\sigma$-weakly dense subalgebra $\mathcal A \subset L^\infty(\G)$.  So, fixing a basis vector $u_{ij}^\alpha$, we compute 
\begin{eqnarray*}
&&((\kappa_r \circ S_\varphi \circ \kappa_r) \otimes \id_{L^\infty(\G)})\Delta_r(u_{ij}^\alpha) = \sum_{k=1}^{d_\alpha} (\kappa_r \circ S_\varphi \circ \kappa_r)u_{ik}^\alpha \otimes u_{kj}^\alpha \\
&=& \sum_{k=1}^{d_\alpha} (\kappa_r \circ S_\varphi)u_{ki}^{\overline{\alpha}} \otimes u_{kj}^\alpha \ \ \ \ (\textrm{using (\ref{eqn_kappa_action})}) \\
&=& \sum_{k=1}^{d_\alpha} \kappa_r \Big( \sum_{l=1}^{d_\alpha}\varphi(u_{kl}^{\overline{\alpha}})u_{li}^{\overline{\alpha}} \Big) \otimes u_{kj}^\alpha \ \ \ \ (\textrm{using (\ref{eqn_action_S})}) \\
&=& \sum_{1 \le k,l \le d_\alpha} \varphi(u_{kl}^{\overline{\alpha}}) u_{il}^{\alpha} \otimes u_{kj}^\alpha. 
\end{eqnarray*}  Therefore,
\begin{eqnarray*}
Q_\varphi u_{ij}^\alpha &=& \Delta_r^{-1} \circ E((\kappa_r \circ S_\varphi \circ \kappa_r) \otimes \id_{L^\infty(\G)})\Delta_ru_{ij}^\alpha)\\
&=& \Delta_r^{-1}\Big( \sum_{1 \le k,l \le d_\alpha} \varphi(u_{kl}^{\overline{\alpha}}) E(u_{il}^{\alpha} \otimes u_{kj}^\alpha)\Big) \\
&=& \Delta_r^{-1}\Big( \sum_{1 \le k,l \le d_\alpha} \varphi(u_{kl}^{\overline{\alpha}}) \frac{\delta_{l,k}}{d_\alpha}\Delta_r(u_{ij}^\alpha)\Big) \ \ \ \ (\textrm{using Lemma \ref{lem_action_cond_exp}}) \\
&=& \sum_{k=1}^{d_\alpha} \frac{\varphi(u_{kk}^{\overline{\alpha}})}{d_\alpha} u_{ij}^\alpha = \frac{\varphi\Big(\sum_{k=1}^{d_\alpha} u_{kk}^{\overline{\alpha}} \Big)}{d_\alpha} u_{ij}^\alpha = \frac{\varphi(\chi_{\overline{\alpha}})}{d_\alpha} u_{ij}^\alpha \\ 
&=& \frac{\psi(\chi_{\overline{\alpha}})}{d_\alpha} u_{ij}^\alpha \ \ \ \ (\textrm{since $\varphi$ is a Hahn-Banach extension of $\psi$}) \\
&=& T_\psi u_{ij}^\alpha. 
\end{eqnarray*} \hfill $\square$

\begin{rem} As alluded to above, the composition (\ref{eqn_averaging}) can be viewed as a type of average of the convolution operator $S_\varphi$, with respect to the Haar trace.  Indeed, consider the classical situation where $\G = (C(G), \Delta)$ is the usual C$^\ast$-bialgebra of continuous functions on a compact group $G$, with coproduct $\Delta$ given by $\Delta f(x,y) = f(xy)$, ($x,y \in G, \ f\in C(G)$).  Let $\mu \in C(G)^* = M(G)$  be a complex regular Borel measure, and for each $g \in G$, let $\mu^g \in M(G)$ be the measure defined by 
\begin{eqnarray*}
\int_{G} f(x) \textrm{d}\mu^g(x) := \int_{G} f(g^{-1}x^{-1}g)\textrm{d} \mu(x) && (f \in C(G)).
\end{eqnarray*}  Then a simple calculation shows that (\ref{eqn_averaging}) applied to $S_\mu$ corresponds to the averaging map $$S_\mu \mapsto Q_\mu= \int_GS_{\mu^g} \textrm{d}g,$$ where $\textrm{d}g$ denotes the Haar probability measure on $G$.  Variants of the type of composition given in (\ref{eqn_averaging}) were first used in the context of discrete group C$^\ast$-algebras by Haagerup \cite{Ha0}, and in the setting of discrete Kac algebras by Kraus and Ruan \cite{KrRu}.  
\end{rem}

\section{The Haagerup Approximation Property for $O_N^+$ and $U_N^+$} \label{section_orthogonal}

In this section, we prove that the free orthogonal quantum groups $\{O_N^+\}_{N \ge 2}$ and the free unitary quantum groups $\{U_N^+\}_{N\ge 2}$  have the HAP.  We will first treat the orthogonal case, and then use this result together with some free product constructions to deal with the unitary case.  

\subsection{The Orthogonal Case}  

The corepresentation theory of the quantum group $O_N^+$ was first studied by Banica \cite{Ba0}.  It turns out to be quite similar to the representation theory of the compact group $SU(2)$.  In the following theorem we have collected the results that we will need from \cite{Ba0}.  

\begin{thm} (\cite{Ba0})  \label{thm_banica_orthogonal} For any fixed $N \ge 2$, there is a maximal family irreducible corepresentations of $O_N^+$, labeled by the non-negative integers, say $$\{V^n = [v^n_{ij}]_{1 \le i,j \le d_n^{(N)}}\}_{n=0}^\infty,$$ which has the following properties: 
\begin{enumerate}
\item $V^0 = 1_{A_o(N)}$ is the trivial corepresentation of $O_N^+$, and $V^1 = V$ is the fundamental corepresentation of $O_N^+$.
\item $\overline{V^n} \cong V^n$ for all $n \ge 0$.
\item The family $\{V^n\}_{n \ge 0}$ satisfies the fusion rules \begin{eqnarray} \label{eqn_fusion_rules} 
V^r \boxtimes V^s \cong \bigoplus_{l = 0}^{\min\{r,s\}} V^{r+s - 2l} && (r,s \in \N \cup \{0\}).
\end{eqnarray}
\item For $N=2$, $d_n^{(2)} = n+1$, and for $N \ge 3$, 
\begin{eqnarray} \label{eqn_dim_reps} 
d_n^{(N)} = \frac{q(N)^{n+1}-q(N)^{-n-1}}{q(N)-q(N)^{-1}},
\end{eqnarray}
where $q(N)$ is defined by $q(N)+ q(N)^{-1} = N$.
\item The irreducible characters $\{\chi_n = (Tr \otimes \id_{A_o(N)})V^n\}_{n \ge 0} \subset A_o(N)$ are self-adjoint, and satisfy the recursion relations \begin{eqnarray} \label{eqn_recursion_characters} \chi_1 \chi_n = \chi_{n+1} + \chi_{n-1} && (n \ge 1).
\end{eqnarray} In particular, we have the equality of $\ast$-algebras \begin{eqnarray} \label{eqn_equality_star_algs}
\textrm{alg}\langle \chi_n: n \ge 0 \rangle = \textrm{alg}\langle 1, \chi_1 \rangle.
\end{eqnarray}
\item \label{spec} The spectral measure $\mu_{\chi_1}$ of $\chi_1$ relative to the Haar state $h$ on $A_o(N)$ is given by Wigner's semicircular law:  
\begin{eqnarray*}
\text{d}\mu_{{\chi_1}}(t) = 1_{[-2,2]}(t) \frac{\sqrt{4 -t^2}}{2 \pi} \text{d}t.
\end{eqnarray*} 
\end{enumerate} 
\end{thm}
Observe that for each $N \ge 2$, part (\ref{spec}) of the above theorem shows that the spectrum $\sigma_{C(O_N^+)} (\pi_{h}(\chi_1))$ of $\pi_{h}(\chi_1)$ in the reduced C$^\ast$-algebra $C(O_N^+)$ is the interval $[-2,2]$.   In the next lemma, we determine the spectrum of the operator $\chi_1$ in the \textit{full} C$^\ast$-algebra $A_o(N)$. 

\begin{lem} \label{lem_spec_fund_char_full_algebra}
Let $N \ge 2$.  Then the spectrum $\sigma_{A_o(N)}(\chi_1)$ of $\chi_1$ in the full C$^\ast$-algebra $A_o(N)$ is the interval $[-N,N]$.   
\end{lem}
This result is perhaps well known.  Since we we could not find a reference, we include a proof.

\begin{proof} The fundamental corepresentation $V^1 = V = [v_{ij}]_{1 \le i,j \le N} \in M_N(A_o(N))$ is an $N \times N$ orthogonal matrix. Therefore $$\|\chi_1\|_{A_o(N)} = \Big\|\sum_{i =1}^Nv_{ii}\Big\|_{A_o(N)} \le \sum_{i=1}^N \|v_{ii}\|_{A_o(N)} \le N,$$  and it follows that $\sigma_{A_o(N)}(\chi_1) \subseteq [-N,N]$.  On the other hand, the universal property of $A_o(N)$ implies that there exists a surjective C$^\ast$-homomorphism \begin{eqnarray*} \pi:A_o(N) \to C(O_N), && \pi(v_{ij}) = o_{ij},\end{eqnarray*} where $C(O_N)$ is the C$^\ast$-algebra of continuous functions on the orthogonal group $O_N$, and $\{o_{ij}\}_{1 \le i,j \le N}$ are the matrix coordinate functions on $O_N$.  In particular, $\pi(\chi_1)$ is just the fundamental character of $O_N$, which is well known to have spectrum equal to $[-N,N]$.  Since the spectrum of any element of a C$^\ast$-algebra always contains the spectrum of its image under any $\ast$-homomorphism, we deduce that \begin{eqnarray*} [-N,N] \subseteq \sigma_{A_o(N)}(\chi_1).
\end{eqnarray*} 
\end{proof}

Let $\{u_n\}_{n=0}^\infty$ denote the sequence of (dilated) Chebyshev polynomials of the second kind, which are uniquely determined by the initial conditions \begin{eqnarray} \label{eqn_ic_chebyshev}
u_0(x) = 1, \ u_1(x) = x
\end{eqnarray} and the recursion 
\begin{eqnarray} \label{eqn_recursion_chebyshev}
xu_n(x) = u_{n+1}(x) + u_{n-1}(x) && (n \ge 1).
\end{eqnarray}
From Theorem \ref{thm_banica_orthogonal} and Lemma \ref{lem_spec_fund_char_full_algebra}, we immediately obtain the following.

\begin{cor} \label{cor_C_alg_fund_character}
Let $\mathcal B_N = C^*\langle \chi_n: n \ge 0 \rangle \subset A_o(N)$ denote the unital C$^\ast$-algebra generated by the irreducible characters of $O_N^+$.  Then there is a C$^\ast$-isomorphism $\mc B_N \cong C([-N,N])$ defined by  \begin{eqnarray*}
\chi_n \mapsto u_n|_{[-N,N]} && (n \ge 0). 
\end{eqnarray*}
\end{cor}

\begin{proof}
From (\ref{eqn_equality_star_algs}), we have $\mc B_N = C^*\langle 1, \chi_1 \rangle.$ On the other hand, Lemma \ref{lem_spec_fund_char_full_algebra} and the Gelfand theorem imply that there is a $\ast$-isomorphism $\pi: C^*\langle 1, \chi_1 \rangle \to C([-N,N])$ determined by $\pi(1_{A_o(N)}) = 1|_{[-N,N]} = u_0|_{[-N,N]}$ and $\pi(\chi_1) = x|_{[-N,N]} = u_1|_{[-N,N]}$.      

Applying $\pi$ to the recursion relation (\ref{eqn_recursion_characters}) therefore yields 
\begin{eqnarray*}
x[\pi(\chi_n)](x) = [\pi(\chi_{n+1})](x)  + [\pi(\chi_{n-1})](x)  && (x \in [-N,N], \ n \ge 1).
\end{eqnarray*}
Comparing with the initial conditions (\ref{eqn_ic_chebyshev}) and the recursion (\ref{eqn_recursion_chebyshev}) for the Chebyshev polynomials of the second kind, we conclude that $\pi(\chi_{n}) = u_n|_{[-N,N]}$ for all $n \ge 0$.
\end{proof}

We now proceed toward proving that $O_N^+$ has the HAP for all $N \ge 2$.  Note that when $N=2$, $O_2^+$ is a co-amenable compact quantum group \cite[Corollaire 1]{Ba0}, and it follows from \cite[Theorem 1.1]{BeMuTu} that $C(O_2^+)$ is a nuclear C$^\ast$-algebra and $L^\infty(O_2^+)$ is an injective von Neumann algebra.  Since injectivity implies the HAP for a finite von Neumann algebra (see \cite{Jo1}), there is nothing to prove when $N=2$.  Assume for the remainder of the section that $N \ge 3$ is fixed.  As in Section \ref{section_general_results}, we will denote by $\{p_n\}_{n \ge 0}$ the orthogonal family of projections $p_n:L^2(O_N^+) \to L^2_n(O_N^+) = \textrm{span}\{u_{ij}^n: 1 \le i,j \le d_n^{(N)}\}$. 

\begin{prop} \label{prop_net_of_maps}
Fix a number $2 < t_0 <3$, and for each $t \in [t_0, N]$ consider the sequence $$\Big\{ \frac{u_n(t)}{u_n(N)} \Big\}_{n = 0}^{\infty}.$$  
\begin{enumerate}
\item \label{decay} There exists a constant $C_{t_0} > 0$ (only depending on $t_0$) such that 
\begin{eqnarray} \label{eqn_coeff_ineq}
0 < \frac{u_n(t)}{u_n(N)} \le C_{t_0}\Big(\frac{t}{N}\Big)^n && (t \in [t_0, N], \ n \ge 0).
\end{eqnarray} 
\item \label{compact} For each $t \in [t_0,N)$, the operator $$\hat{T}_t = \sum_{n \ge 0} \frac{u_n(t)}{u_{n}(N)}p_n \in B(L^2(O_N^+)),$$ is compact.
\item \label{cp}  For each $t \in [t_0,N]$, $\hat{T}_t$ restricts to a NUCP, $h$-preserving map $T_t \in CB(L^\infty(O_N^+))$ with the property that $T_t(C(O_N^+)) \subseteq C(O_N^+)$.
\end{enumerate}   
\end{prop}

\begin{proof}
(\ref{decay}):  For $t > 2$, let \begin{eqnarray} \label{eqn_q} q(t) = \frac{t + \sqrt{t^2-4}}{2}. \end{eqnarray}  Then then $q$ is an increasing function of $t$, and $q(t) + q(t)^{-1} = t$.  Using induction and the recursion (\ref{eqn_recursion_chebyshev}), it follows that for all $t > 2$ and $n \ge 1$,
\begin{eqnarray} \label{eqn_q_cheby}
0 < u_n(t) = u_n(q(t)+q(t)^{-1}) = \frac{q(t)^{n+1}-q(t)^{-n-1}}{q(t)-q(t)^{-1}}.
\end{eqnarray}
In particular, for any $t \ge t_0$ and any $n \ge 1$, we have 
\begin{eqnarray*}
0 &<& \frac{u_n(t) }{u_n(N)}  = \Bigg(\frac{q(t)}{q(N)}\Bigg)^n \Bigg(\frac{(1-q(t)^{-2n-2})(1 -q(N)^{-2})}{(1-q(t)^{-2})(1-q(N)^{-2n-2})}\Bigg) \\
&\le& \Bigg(\frac{q(t)}{q(N)}\Bigg)^n \Bigg(\frac{1 -q(N)^{-2}}{(1-q(t)^{-2})(1-q(N)^{-2n-2})}\Bigg) \\
&\le& \Bigg(\frac{q(t)}{q(N)}\Bigg)^n(1-q(t)^{-2})^{-1} \\
&& (\textrm{since $1 - q(N)^{-2} \le 1-q(N)^{-2n-2}$})\\
&\le&  \Bigg(\frac{q(t)}{q(N)}\Bigg)^n(1-q(t_0)^{-2})^{-1} \\
&& (\textrm{since $q$ is increasing}) \\
&=& (1-q(t_0)^{-2})^{-1} \Bigg(\frac{1 + \sqrt{1 - 4/t^2}}{1 + \sqrt{1-4/N^2}}\Bigg)^n \Bigg( \frac{t}{N}\Bigg)^n \\
&\le&  (1-q(t_0)^{-2})^{-1}\Bigg( \frac{t}{N}\Bigg)^n.
\end{eqnarray*}
Taking $C_{t_0} = (1-q(t_0)^{-2})^{-1}$, we obtain (\ref{eqn_coeff_ineq}). 

\noindent (\ref{compact}): Part (\ref{decay}) in particular shows that the sequence $\Big\{ \frac{u_n(t)}{u_n(N)} \Big\}_{n = 0}^{\infty}$ belongs to $c_0(\N \cup \{0\})$ for each $t \in [t_0,N)$.  Since the projections $p_n$ $(n \ge 0)$ are all finite rank and mutually orthogonal, the sum $$\hat{T}_t = \sum_{n \ge 0} \frac{u_n(t)}{u_{n}(N)}p_n \in B(L^2(O_N^+)),$$ is compact. 

\noindent (\ref{cp}):  For each $t \in [-N,N]$, let $\psi_t \in C([-N,N])^*$ denote the state given by point evaluation at $t$.  Using the isomorphism $$C([-N,N]) \cong \mc B_N = C^*\langle \chi_n: n \ge 0 \rangle,$$ given by Corollary \ref{cor_C_alg_fund_character}, we may regard $\psi_t$ as the state in $\mc B_N^*$ given by \begin{eqnarray*} \psi_t(\chi_n) = u_n(t) && (n \ge 0). \end{eqnarray*}
Theorem \ref{thm_central_cp_maps} therefore implies that the map $\hat{T}_{\psi_t} = \sum_{n \ge 0} \frac{u_n(t)}{d_n^{(N)}}p_n \in B(L^2(O_N^+))$ restricts to a NUCP $h$-preserving map $T_{\psi_t} \in CB(L^\infty(O_N^+)),$ and that $T_{\psi_t}$ preserves $C(O_N^+)$.  Since $d_n^{(N)} = u_n(N)$ (cf. equations (\ref{eqn_dim_reps}) and (\ref{eqn_q_cheby})), we have $\hat{T}_t = \hat{T}_{\psi_t}$ for all $t \in [t_0,N]$, completing the proof.
\end{proof}

It is now easy to deduce the HAP for $O_N^+$.

\begin{thm} \label{thm_HAP_orth}
For each $N \ge 2$, $O_N^+$ has the HAP.
\end{thm}

\begin{proof}
As discussed before, we may assume that $N \ge 3$ is fixed.  Let $\{T_t\}_{t \in [t_0, N)} \subset CB(L^\infty(O_N^+))$ be the net of NUCP $h$-preserving maps constructed in Proposition \ref{prop_net_of_maps} above.  Then, for each $t \in [t_0,N)$, the $L^2$-extension $\hat{T}_t \in B(L^2(O_N^+))$ is compact.  It remains to show that 
\begin{eqnarray} \label{eqn_limit}
\lim_{t \to N}\|\hat{T}_t\hat{x} - \hat{x}\|_{L^2(O_N^+)} = 0 && (x \in L^\infty(O_N^+)).
\end{eqnarray} 
Since the family $\{\hat{T}_t\}_{t <N} \subset B(L^2(O_N^+))$ is uniformly bounded (by $1$), it suffices by linearity and $L^2$-density to prove (\ref{eqn_limit}) for all matrix elements $x = v_{ij}^n$, ($n \ge 0$, $1 \le i,j \le d_n^{(N)}$).  But in this case, we have $$\lim_{t \to N} \|\hat{T}_t v_{ij}^n - v_{ij}^n\|_{L^2(O_N^+)} = \lim_{t \to N} \Big(\frac{u_n(t)}{u_n(N)} - 1\Big) (d_n^{(N)})^{-1/2} = 0.$$  
\end{proof}

\begin{rem}\label{rem_pt_nm_convergence_orth}
An easy modification of the above proof (replacing $L^2$-norms with $C(O_N^+)$-norms) also shows that $T_t|_{C(O_N^+)} \to \id_{C(O_N^+)}$, in the point-norm topology of $B(C(O_N^+))$.
\end{rem}

\subsection{The Unitary Case} 

Let us now consider the free unitary quantum groups $\{U_N^+\}_{N \ge 2}$.  It turns out that using free product constructions, the HAP for $U_N^+$ is really a consequence of the HAP for $O_N^+$.  

Let $A$ and $B$ be unital C$^\ast$-algebras, we will denote their C$^\ast$-algebraic free product by $A*B$ (\cite[Chapter 25]{Pi}).  If $\varphi$ and $\psi$ are faithful states on $A$ and $B$, respectively, we will denote by $A*_{\textrm{red}}B$ the reduced free product of $A$ and $B$, taken with respect to the free product state $\varphi*\psi$ (see \cite[Lecture 7]{NiSp}).  If $A$ and $B$ are von Neumann algebras with normal faithful states $\varphi$ and $\psi$, respectively, then we always take $A*_{\textrm{red}}B$ to be the Neumann algebraic reduced free product \cite{BlDy}.  We refer to \cite{Ch} and \cite{BlDy} for the properties of reduced free products of (normal) completely positive maps.  

We begin with the notion of free complexification for compact matrix quantum groups, which was introduced in \cite{Ba2}.  Recall that a \textit{compact matrix quantum group} $\G = (A, U)$ is simply a CQG $(A,\Delta)$, where $U = [u_{ij}]_{1 \le i,j \le N} \in M_N(A)$ is a distinguished \textit{fundamental} unitary corepresentation, with the property that its matrix elements generate $A$ as a C$^\ast$-algebra \cite{Wo}. 

\begin{defn} \label{defn_free_complexification}
Let $\G = (A, U)$ be a compact matrix quantum group, with fundamental corepresentation $U=[u_{ij}]_{1 \le i,j \le N} \in M_N(A)$.  Let $\T$ denote the unit circle in $\C$, and let $z = \id_\T \in C(\T)$ denote the canonical unitary generator of $C(\T)$.  The \textit{free complexification} of $\G$ is the compact matrix quantum group $\tilde{\G} = (\tilde{A}, \tilde{U})$ where $\tilde{A} \subset C(\T)*A$ is the unital C$^*$-algebra generated by the family $\{ zu_{ij}: 1 \le i,j \le N\} \subset C(\T)*A,$ and $\tilde{U} = [zu_{ij}]_{1 \le i,j \le N} \in M_N(\tilde{A})$.  The fact that $\tilde{\G}$ is actually a compact matrix quantum group is proved in \cite{Ba2}.
\end{defn}

Denote by $\tau$ the Haar state on $L^\infty(\T)$, and $h_\G$ the Haar state on $L^\infty(\G)$.  In \cite{Ba2} it is shown that $L^\infty(\tilde{\G})$ (respectively $C(\tilde{\G})$) is the von Neumann subalgebra of $L^\infty(\T)*_{\textrm{red}} L^\infty(\G)$ (respectively C$^\ast$-subalgebra of $C(\T)*_{\textrm{red}} C(\G)$) generated by $\{z\pi_{h_\G}(u_{ij}): 1 \le i,j \le N\}$.  The Haar state $h_{\tilde{\G}}$ on $L^\infty(\tilde{\G})$ is just the restriction of the free product state $\tau*h_\G$ to this von Neumann subalgebra.

Here is the main result of this section.

\begin{thm} \label{thm_HAP_free_complexification}
Let $\G$ be a compact matrix quantum group of Kac type.  If $\G$ has the HAP, then so does its free complexification $\tilde{\G}$.
\end{thm}    

\begin{proof}
Suppose $(M_1, \tau_1)$ and $(M_2,\tau_2)$ are finite von Neumann algebras, each with the HAP.  Then the reduced free product $(M_1*_{\textrm{red}}M_2, \tau_1*\tau_2)$ also has the HAP \cite[Proposition 3.9]{Boc}.  The main idea in \cite{Boc} is that if $\{\Phi^{(i)}_{\lambda_i}:M_i \to M_i\}_{\lambda_i \in \Lambda_i}$ is a net giving the HAP for $(M_i,\tau_i)$, then the reduced free product net $$\{\Phi^{(1)}_{\lambda_1}*_{\textrm{red}}\Phi^{(2)}_{\lambda_2} :M_1*_{\textrm{red}}M_2 \to M_1*_{\textrm{red}}M_2\}_{(\lambda_1, \lambda_2) \in \Lambda_1 \times \Lambda_2},$$ gives the HAP for $(M_1*_{\textrm{red}}M_2, \tau_1*\tau_2)$.

In particular, $(L^\infty(\T)*_{\textrm{red}}L^\infty(\G), \tau*h_\G)$ has the HAP.  Since $L^\infty(\tilde{\G})$ is a von Neumann subalgebra of $(L^\infty(\T)*_{\textrm{red}}L^\infty(\G), \tau*h_\G)$, and $\tau*h_\G$ is a faithful normal trace, there exists a unique normal, faithful $\tau*h_\G$-preserving conditional expectation $E: L^\infty(\T)*_{\textrm{red}}L^\infty(\G) \to L^\infty(\tilde{\G})$.  Thus, if $\{\Phi_{\lambda}\}_{\lambda \in \Lambda}$ is any net giving the HAP for $(L^\infty(\T)*_{\textrm{red}}L^\infty(\G), \tau*h_\G)$, then the net $\{E \circ \Phi_{\lambda}|_{L^\infty(\tilde{\G})}\}_{\lambda \in \Lambda}$ gives the HAP for $(L^\infty(\tilde{\G}),h_{\tilde{\G}})$.  
\end{proof}

As a specific instance of Theorem \ref{thm_HAP_free_complexification}, we deduce the HAP for the free unitary quantum groups. 

\begin{cor} \label{cor_HAP_unitary}
For each $N \ge 2$, $U_N^+$ has the HAP.
\end{cor}

\begin{proof}
In \cite[Th\'eor\`eme 1]{Ba} (see \cite[Theorem 9.2]{BaCo} for another proof), it is shown that $U_N^+ = \tilde{O}_N^+$ for all $N \ge 2$.  Therefore it follows from Theorem \ref{thm_HAP_orth} and Theorem \ref{thm_HAP_free_complexification} that $U_N^+$ has the HAP.   
\end{proof}

\begin{rem} \label{rem_explicit_net_unitary_case}
Of course, Corollary \ref{cor_HAP_unitary} does not produce an explicit net of NUCP $h$-preserving maps on $L^\infty(U_N^+)$ yielding the HAP for $U_N^+$.  Since we will need one in Section \ref{section_MAP_uni}, let us now produce one.  For each $r \in [0,1)$, let $p_r(e^{i\theta}) = \sum_{n \in \Z} r^{|n|}e^{in\theta} \in C(\T)$ denote the Poisson kernel.  It is well known that $p_r > 0$ and $\frac{1}{2\pi}\int_{-\pi}^\pi p_r(e^{i\theta}) \text{d}\theta =1$ for all $r \in [0,1)$, and therefore convolution operator \begin{eqnarray*}
f \mapsto P_r f := p_r*f && (f \in L^\infty(\T)),
\end{eqnarray*} defines a NUCP Haar trace-preserving map on $L^\infty(\T)$, such that $P_r( C(\T)) \subseteq C(\T)$.  Furthermore, if $z = \id_\T$ denotes the canonical generator of $L^\infty(\T)$, then $\{z^n\}_{n \in \Z}$ is an orthonormal basis for $L^2(\T)$ on which $P_r$ is given by $P_r(z^n) = r^{|n|}z^n$.  From this, it is clear that $P_r$ extends to a compact (multiplication) operator on $L^2(\T)$, and $\lim_{r \to 1}P_r f = f$ for all $f \in L^2(\T)$.  Thus $\{P_r\}_{r \in [0,1)}$ is a net of NUCP trace-preserving maps giving the HAP for $\T$.

Now let $\{\Phi_{r,t}\}_{r \in[0,1), t \in [t_0,N)} \subset CB(L^\infty(\T)*_{\textrm{red}}L^\infty(O_N^+))$ be the reduced free product net given by \begin{eqnarray} \label{eqn_Phi_r,t}
\Phi_{r,t} = P_r*_{\textrm{red}}T_t,
\end{eqnarray}
where $T_t \in CB(L^\infty(O_N^+))$ is the NUCP $h$-preserving map constructed in Proposition \ref{prop_net_of_maps}.  Then the net $\{\Phi_{r,t}\}_{r \in[0,1), t \in [t_0,N)}$ gives the HAP for $(L^\infty(\T)*_{\textrm{red}}L^\infty(O_N^+), \tau*h)$.  In Proposition \ref{prop_restriction}, we will see that in fact $\Phi_{r,t}(L^\infty(U_N^+)) \subseteq L^\infty(U_N^+)$ (and also $\Phi_{r,t}(C(U_N^+)) \subseteq C(U_N^+)$), so $\{\Phi_{r,t}|_{L^\infty(U_N^+)}\}_{r \in[0,1), t \in [t_0,N)}$ also yields the HAP for $U_N^+$.

In our proof of the metric approximation property for $C(U_N^+)$ in Section \ref{section_MAP_uni}, we will make extensive use of the maps $\Phi_{r,t}$ defined above. 
\end{rem}

\section{The Reduced C$^\ast$-Algebras and the Metric Approximation Property} \label{section_MAP}

In this section, we use the HAP for $O_N^+$ and $U_N^+$ to show that the reduced C$^\ast$-algebras $C(O_N^+)$ and $C(U_N^+)$ have the metric approximation property.  Let us begin by recalling this notion. 

\begin{defn}\label{def_MAP}
A Banach space $Y$ has the \textit{metric approximation property (MAP)} if there exists a net $\{\Phi_\lambda\}_{\lambda \in \Lambda} \subset B(Y)$ of finite rank contractions converging to the identity map in the point-norm topology of $B(Y)$.  That is, for all $y \in Y$ $$\lim_{\lambda \in \Lambda}\|\Phi_\lambda y - y\|_Y = 0.$$ \end{defn}     

We will treat the orthogonal case and the unitary case separately, starting with the orthogonal case. 

\subsection{The Metric Approximation Property for $C(O_N^+)$}   Since $C(O_2^+)$ is a nuclear C$^\ast$-algebra (and therefore trivially has the MAP), we will assume for the remainder of this section that $N \ge 3$. 

In Proposition \ref{prop_net_of_maps}, we constructed a net of UCP maps $\{T_t\}_{t \in [t_0,N)} \subset CB(C(O_N^+))$, each of which is $L^2$-compact, with the property that $\lim_{t \to N} T_t = \id_{C(O_N^+)}$ in the point-norm topology (see Remark \ref{rem_pt_nm_convergence_orth}).  Furthermore, Proposition \ref{prop_net_of_maps}(\ref{decay}) shows that for fixed $t$, the operator $T_t$ has rapid (i.e. exponential) decay properties as the corepresentation index $n$ tends to $\infty$.  To obtain an appropriate net of finite rank contractions on $C(O_N^+)$ converging to $\id_{C(O_N^+)}$ in the point-norm topology, the idea is to take certain finite rank truncations of the maps $T_t$.  It turns out that we can meaningfully control the norms of these truncated maps using: $(1)$ the explicit decay rates of $T_t$ given in Proposition \ref{prop_net_of_maps}(\ref{decay}), and $(2)$ the fact that the quantum groups $O_N^+$ ($N \ge 3$) have the \textit{property of rapid decay} (property RD).  The property of rapid decay, defined for CQGs by Vergnioux \cite{Ve}, amounts to the existence of a Haagerup-type inequality comparing the $L^2$- and $L^\infty$-norms on $L^\infty(\G)$.  This approach to proving the MAP dates back to Haagerup \cite{Ha} in the setting of free groups.   We state Vergnioux's Haagerup inequality for $O_N^+$ in the following theorem.

\begin{thm} \label{thm_property_RD_orth} (\cite[Section 4]{Ve})
For each $N \ge 3$, there is a constant $D_N > 0$ such that for any $n \ge 0$ and any $x \in L^2_n(O_N^+)$, we have 
\begin{eqnarray*}
\|x\|_{L^2(O_N^+)} \le \|x\|_{L^\infty(O_N^+)} \le D_N (n+1)\|x\|_{L^2(O_N^+)}.
\end{eqnarray*}
\end{thm}       
In other words, Theorem \ref{thm_property_RD_orth} says that on the subspaces $L^2_n(O_N^+) = \textrm{span}\{v_{ij}^n: 1 \le i,j \le d_n^{(N)}\} \subset L^2(O_N^+)$, the $L^\infty$-norm grows at most linearly in $n$, relative to the $L^2$-norm.

An immediate consequence of this Haagerup inequality for $O_N^+$ is the following criterion for determining when certain linear maps $T \in B(L^2(O_N^+))$ are \textit{ultracontractive} (i.e., map $L^2(O_N^+)$ boundedly into $C(O_N^+)$.)

\begin{prop} \label{prop_ultracontractive}
Let $\{a_n\}_{n \ge 0} \subset \C$ be a bounded sequence, and consider the operator $T = \sum_{n \ge 0} a_n p_n \in B(L^2(O_N^+))$.  If $$k_a:= \sup_{n \ge 0} (n+1)^2 |a_n| < \infty,$$
then $T(L^2(O_N^+)) \subseteq C(O_N^+)$ and $\|T\|_{L^2 \to L^\infty} \le \frac{\pi D_N k_a}{\sqrt{6}}$, where $D_N>0$ is the constant given in Theorem \ref{thm_property_RD_orth}.  
\end{prop} 

\begin{proof}
Let $x \in L^2(O_N^+)$.  Then we have 
\begin{eqnarray*}
&& \|Tx\|_{L^\infty(O_N^+)} = \Big\| \sum_{n \ge 0}a_np_nx \Big\|_{L^\infty(O_N^+)} \\
&\le& \sum_{n \ge 0}|a_n| \cdot \|p_nx\|_{L^\infty(O_N^+)} \\
&\le& \sum_{n \ge 0} |a_n|D_N(n+1) \|p_nx\|_{L^2(O_N^+)} \\
&\le & \sum_{n \ge 0} \frac{k_a D_N}{n+1}\|p_nx\|_{L^2(O_N^+)} \\
&\le& k_a D_N \Big( \sum_{n \ge 0} \frac{1}{(n+1)^2}\Big)^{1/2}\Big( \sum_{n \ge 0} \|p_nx\|_{L^2(O_N^+)}^2\Big)^{1/2} \\
&=& \frac{\pi D_N k_a}{\sqrt{6}} \|x\|_{L^2(O_N^+)},
\end{eqnarray*}
where in the second inequality we have used property RD, and in the last inequality we used the Cauchy-Schwarz inequality.  Therefore $T(L^2) \subseteq L^\infty$ and $\|T\|_{L^2 \to L^\infty} \le \frac{\pi D_N k_a}{\sqrt{6}}$.  To show that $T(L^2(O_N^+)) \subseteq C(O_N^+)$, note that $T \mc A \subseteq \mc A$, where $\mc A$ is the usual dense Hopf $\ast$-subalgebra of $C(O_N^+)$.  Since $\mc A$ is also dense in $L^2(O_N^+)$ and $T:L^2 \to L^\infty$ is continuous, $T(L^2(O_N^+)) \subset \overline{\mc A}^{\|\cdot\|_{L^\infty(O_N^+)}} = C(O_N^+)$. 
\end{proof}

The MAP now follows.

\begin{thm} \label{thm_MAP_orth}
For each $N \ge 3$, $C(O_N^+)$ has the MAP. 
\end{thm}

\begin{proof}
Fix $N \ge 3$, and let $\{T_t\}_{t \in [t_0,N)} \subset CB(C(O_N^+))$ be the family of UCP maps constructed in Proposition \ref{prop_net_of_maps}.  For each $m \in \N$, let $T_{t,m}:C(O_N^+) \to C(O_N^+)$ be the $m$th-order truncation of $T_t$ given by 
\begin{eqnarray*}
T_{t,m} = \sum_{n = 0}^m \frac{u_n(t)}{u_n(N)} p_n.
\end{eqnarray*}  
Since each member of the family $\{p_n\}_{n \ge 0}$ is a finite rank projection, it follows that the finite sums $\{T_{t, m}\}_{m \ge 0, t \in [t_0,N)}$ are also finite rank.  Furthermore, for fixed $t \in [t_0,N)$ we have 
\begin{eqnarray}
\|T_t - T_{t,m}\|_{B(C(O_N^+))} &=& \Big\| \sum_{n \ge m+1}\frac{u_n(t)}{u_n(N)} p_n\Big\|_{B(C(O_N^+))} \\
\label{inclusion}&\le& \Big\| \sum_{n \ge m+1}\frac{u_n(t)}{u_n(N)} p_n\Big\|_{B(L^2(O_N^+),C(O_N^+))} \\
\label{prop}&\le&  \frac{\pi D_N}{\sqrt{6}} \sup_{n \ge m+1} (n+1)^2 \frac{u_n(t)}{u_n(N)} \\
\label{propp}&\le& \frac{\pi D_N}{\sqrt{6}} \sup_{n \ge m+1} (n+1)^2 C_{t_0} \Big(\frac{t}{N}\Big)^n,
\end{eqnarray}
where in (\ref{inclusion}) we have used the fact that the inclusion $C(O_N^+) \subseteq L^2(O_N^+)$ is a contraction, in (\ref{prop}) we have used Proposition \ref{prop_ultracontractive}, and in (\ref{propp}) we have used Proposition \ref{prop_net_of_maps}(\ref{decay}).  From the last inequality, it follows that \begin{eqnarray} \label{n1}
\lim_{m \to \infty}\|T_t - T_{t,m}\|_{B(C(O_N^+))} =0 && (t \in [t_0,N)),  
\end{eqnarray}
and in particular
\begin{eqnarray} \label{n2}
\lim_{m \to \infty} \|T_{t,m}\|_{B(C(O_N^+))} = \|T_{t}\|_{B(C(O_N^+))} = 1 && (t \in [t_0,N)).
\end{eqnarray}

Put $\tilde{T}_{t,m} = \|T_{t,m}\|_{B(C(O_N^+))}^{-1} T_{t,m}$.  Then $\{\tilde{T}_{t,m}\}_{t \in [t_0,N), m \in \N} \subset B(C(O_N^+))$ is a family of finite rank contractions, and by (\ref{n1})--(\ref{n2}), 

\begin{eqnarray} \label{n3}
\lim_{m \to \infty}\|T_t - \tilde{T}_{t,m}\|_{B(C(O_N^+))} =0 && (t \in [t_0,N)).  
\end{eqnarray} 
We now claim that the identity map $\id_{C(O_N^+)}$ is contained in the point-norm closure closure of $\{\tilde{T}_{t,m}\}_{t \in [t_0,N), m \in \N} \subset B(C(O_N^+))$.  This will complete the proof, since we can then extract a net of finite rank contractions (in fact a sequence) from $\{\tilde{T}_{t,m}\}_{t \in [t_0,N), m \in \N}$ with the required property.  To prove this last claim, first note that by (\ref{n3}), $\{T_t\}_{t \in [t_0,N)}$ is contained in the norm ($\Rightarrow$ point-norm) closure of $\{\tilde{T}_{t,m}\}_{t \in [t_0,N), m \in \N}$.  On the other hand, since $\lim_{t \to N}\frac{u_n(t)}{u_n(N)} = 1$ for all $n \ge 0$, we have $\lim_{t \to N}\|T_t x - x\|_{C(O_N^+)} = 0$ for all $x \in \mc A$, the dense Hopf $\ast$-subalgebra of $C(O_N^+)$.  Since $\|T_t\|_{B(C(O_N^+))} = 1$ for all $t \in [t_0,N)$, the density of $\mc A \subseteq C(O_N^+)$ implies that this limit remains valid for all $x \in C(O_N^+)$.   Therefore $$\id_{C(O_N^+)} \in \overline{\{T_t\}_{t \in [t_0,N)}}^{\textrm{ point-norm}} \subset \overline{\{\tilde{T}_{t,m}\}_{t \in [t_0,N), m \in \N}}^{\textrm{ point-norm}}.$$  
\end{proof}

\subsection{The Metric Approximation Property for $C(U_N^+)$} \label{section_MAP_uni}
To prove the MAP for the C$^\ast$-algebras $C(U_N^+),$ the idea is the same as for $C(O_N^+)$:  we take a net of $L^2$-compact, UCP maps on $C(U_N^+)$ converging to $\id_{C(U_N^+)}$ in the point-norm topology, and consider finite rank truncations of this net.

To begin, we will need to recall some of the facts obtained by Banica \cite{Ba} on the corepresentation theory of $U_N^+$, ($N \ge 2$).  Let $\{V^n\}_{n \ge 0}$ be a complete family of irreducible unitary corepresentations of $O_N^+$ as in Theorem \ref{thm_banica_orthogonal}, and let $z = \id_\T$ denote the canonical unitary generator (i.e. fundamental corepresentation) of $C(\T)$.  Since $U_N^+$ is the free complexification of $O_N^+$ (see Definition \ref{defn_free_complexification}), we may identify $C(U_N^+)$ with the C$^\ast$-subalgebra of $C(\T)*_{\textrm{red}}C(O_N^+)$ generated by the matrix elements of $z \boxtimes V^1 = [zv_{ij}]_{1 \le i,j \le N}$.  

Let $\F_2^+$ denote the free semigroup on two generators $\{g_1, g_2\}$ with unit $e$, and equip $\F_2^+$ with the antimultiplicative involution $g \mapsto \overline{g}$ defined by $$\overline{e} = e, \ \ \ \overline{g_1} = g_2, \ \ \ \overline{g_2} = g_1.$$  Then there exists a maximal family of pairwise inequivalent irreducible unitary corepresentations $\{U^g\}_{g \in \F_2^+}$ of $U_N^+$ such that
\begin{enumerate}
\item $U^e = 1$,
\item $U^{g_1} = z \boxtimes V^1$, the fundamental corepresentation of $U_N^+$,
\item $\overline{U^{g}} \cong U^{\overline{g}}$ for all $g \in \F_2^+$, and 
\item The following fusion rules are satisfied: 
\begin{eqnarray} \label{eqn_fusion_unitary}
U^{g}\boxtimes U^{h} \cong \bigoplus_{\substack{\alpha, \beta, \sigma \in \F_2^+ \\ g = \alpha\sigma, \ h = \overline{\sigma}\beta}} U^{\alpha\beta} && (g,h \in \F_2^+). 
\end{eqnarray}
\end{enumerate} 
For $N \ge 3$, Vergnioux has also shown that $U_N^+$ also has property RD \cite[Section 4]{Ve}.  Let $|\cdot|:\F_2^+ \to \N \cup \{0\}$ denote the word length function on $\F_2^+$ relative to the generating set $\{g_1,g_2\}$.   For each $n \in \N \cup \{0\}$, let $$L^2_n(U_N^+) = \bigoplus_{g \in \F_2^+: |g|=n}L^2_g(U_N^+),$$ where $L^2_g(U_N^+)$ denotes the usual subspace of $L^2(U_N^+)$ spanned by the matrix elements of the corepresentation $U^g$.  Vergnioux's property RD for $U_N^+$ can now be stated as follows.

\begin{thm} \label{thm_property_RD_unit} (\cite[Section 4]{Ve})
For each $N \ge 3$, there is a constant $R_N > 0$ such that for any $n \ge 0$ and any $x \in L^2_n(U_N^+)$, we have 
\begin{eqnarray*}
\|x\|_{L^2(U_N^+)} \le \|x\|_{L^\infty(U_N^+)} \le R_N (n+1)\|x\|_{L^2(U_N^+)}.
\end{eqnarray*}
\end{thm}  

For each $g \in \F_2^+$, let $p_g:L^2(U_N^+) \to L^2_g(U_N^+)$ be the orthogonal projection.  Using Theorem \ref{thm_property_RD_unit}, we have the following ultracontractivity result for $U_N^+$ (which is analogous to Proposition \ref{prop_ultracontractive} for $O_N^+$).  

\begin{prop} \label{prop_ultracontractive_unit}
Fix $N \ge 3$, let $\{a_g\}_{g  \in \F_2^+} \subset \C$ be a bounded family, and consider the operator $T = \sum_{g \in \F_2^+} a_g p_g \in B(L^2(U_N^+))$.  If $$k_a:= \sup_{n \ge 0}\Big\{(n+1)^2 \max_{|g|=n} |a_g|\Big\} < \infty,$$
then $T(L^2(U_N^+)) \subseteq C(U_N^+)$ and $\|T\|_{L^2 \to L^\infty} \le \frac{\pi R_N k_a}{\sqrt{6}}$, where $R_N>0$ is the constant given in Theorem \ref{thm_property_RD_unit}.  
\end{prop}

We will omit the proof of this result since it is almost identical to that of Proposition \ref{prop_ultracontractive}.  

Proposition \ref{prop_ultracontractive_unit} will be essential to our proof of the MAP for $C(U_N^+)$ when $N \ge 3$.  We will treat the case $N =2$ separately at the end of this section.  

\subsubsection{The Case $N \ge 3$.}        
Let $\{\Phi_{r,t}\}_{r \in[0,1), t \in [t_0,N)} \subset CB(L^\infty(\T)*_{\textrm{red}}L^\infty(O_N^+))$ be the family of NUCP maps defined in Remark \ref{rem_explicit_net_unitary_case}.  In the following proposition, we show that each $\Phi_{r,t}$ preserves the subalgebra $L^\infty(U_N^+) \subseteq L^\infty(\T)*_{\textrm{red}}L^\infty(O_N^+)$.

\begin{prop} \label{prop_restriction}
$\Phi_{r,t}(C(U_N^+)) \subseteq C(U_N^+)$ and $\Phi_{r,t}(L^\infty(U_N^+)) \subseteq L^\infty(U_N^+)$. 
\end{prop}

\begin{proof}
For each $g \in \F_2^+$, consider the subspace $L^2_g(U_N^+) \subset C(U_N^+)$ spanned by the matrix elements of the associated irreducible corepresentation $U^g$ of $U_N^+$.  We claim that $\Phi_{r,t}|_{L^2_g(U_N^+)} \in \C \id_{L^2_g(U_N^+)}$.  Assuming that this claim is true, it then follows by linearity, continuity and density, that $\Phi_{r,t}(C(U_N^+)) \subseteq C(U_N^+)$.  Since $\Phi_{r,t}$ is also normal, $\Phi_{r,t}(L^\infty(U_N^+)) \subseteq L^\infty(U_N^+)$.  

We now prove the above claim.  Fix $g \in \F_2^+\backslash \{e\}$ (the case $g=e$ is trivial), and consider the free product compact quantum group $\G = \T*O_N^+$, as defined in \cite{Wa}.  Since $$L^\infty(U_N^+) \subseteq L^\infty(\T)*_{\textrm{red}}L^\infty(O_N^+) =: L^\infty(\G),$$ and the Haar state on $L^\infty(U_N^+)$ is just the restriction of the Haar state on $L^\infty(\G)$, it follows that $U^g$ is also an irreducible corepresentation of $\G$.  Therefore, 
by \cite[Theorem 1.1]{Wa}, there exist numbers $n \in \N$, $l(1), \ldots, l(n) \in \N$, $k(1),k(n+1) \in \Z$, and $k(2), \ldots, k(n) \in \Z \backslash \{0\}$ such that $U^g$ is unitarily equivalent to the alternating tensor product corepresentation $$W = z^{k(1)}\boxtimes V^{l(1)} \boxtimes z^{k(2)} \boxtimes \ldots \boxtimes z^{k(n)} \boxtimes V^{l(n)}\boxtimes z^{k(n+1)} \in M_{d_g}(L^\infty(\G)).$$   

From the above equation, the definition of the reduced free product of NUCP maps \cite{BlDy}, and the definition of $P_r$  and $T_t$ (c.f. Proposition \ref{prop_net_of_maps} and Remark \ref{rem_explicit_net_unitary_case}),  we have
\begin{eqnarray*}
(\id \otimes \Phi_{r,t})W &=& (\id \otimes (P_r*_{\textrm{red}}T_t))W \\
&=& r^{\sum_{s=1}^{n+1}|k(s)|} \prod_{s=1}^n \frac{u_{l(s)}(t)}{u_{l(s)}(N)} W.
\end{eqnarray*}
Since $U^g \cong W$, the same is true for $U^g$: \begin{eqnarray}(\id \otimes \Phi_{r,t})U^g = r^{\sum_{s=1}^{n+1}|k(s)|} \prod_{s=1}^n \frac{u_{l(s)}(t)}{u_{l(s)}(N)} U^g.
\end{eqnarray} That is, $\Phi_{r,t}|_{L^2_g(U_N^+)} = r^{\sum_{s=1}^{n+1}|k(s)|} \prod_{s=1}^n \frac{u_{l(s)}(t)}{u_{l(s)}(N)}\id_{L^2_g(U_N^+)}$.  
\end{proof}

For each $t \in [t_0,N)$, let $q(t) = \frac{t +\sqrt{t^2-4}}{2}$ be the increasing function defined in  (\ref{eqn_q}), let $r(t) = \frac{1-q(t)^{-2}}{1-q(N)^{-2}}$, and consider the UCP map $$\Psi_t := \Phi_{r(t),t}|_{C(U_N^+)} \in CB(C(U_N^+)).$$  From the proof of Proposition \ref{prop_restriction}, we know that $\Psi_t$ acts as a scalar multiple of the identity on each subspace $L^2_g(U_N^+)$, and so its $L^2$-extension $\hat{\Psi}_t$ takes the form 
\begin{eqnarray} \label{eqn_L2_exp}
\hat{\Psi}_t = \sum_{g \in \F_2^+} a_t(g)p_g \in B(L^2(U_N^+)),
\end{eqnarray}
where $\{a_t(g)\}_{g \in \F_2^+} \subset \C$.
In the unitary case, we will use the net $\{\Psi_t\}_{t \in [t_0,N)} \subset CB(C(U_N^+))$ as a replacement for the net $\{T_t\}_{t \in [t_0,N)} \subset CB(C(O_N^+))$ in the orthogonal case.  The following proposition shows that our choice was a good one.

\begin{prop} \label{prop_decay_psi}
The net $\{\Psi_t\}_{t \in [t_0,N)}$ has the following properties: 
\begin{enumerate}
\item $\lim_{t \to N}\Psi_t = \id_{C(U_N^+)}$ in the point-norm topology.
\item For each $g \in \F_2^+$, 
\begin{eqnarray} \label{decay_unitary}
0 < a_t(g) \le C_{t_0} \Big(\frac{t}{N}\Big)^{|g|},
\end{eqnarray}
where $C_{t_0}$ is the constant defined in Proposition \ref{prop_net_of_maps}(\ref{decay}). 
\end{enumerate}
\end{prop}

\begin{proof} Since each map $\Psi_t$ is UCP, it suffices by continuity and density to show that $\lim_{t \to N} \Psi_t x = x$ for each $x \in \mc A \subset C(U_N^+)$, the dense Hopf $\ast$-subalgebra.  By linearity, this reduces to showing that \begin{eqnarray} \label{pointwise_lim} \lim_{t \to N} a_t(g) = 1 && \textrm{pointwise.}
\end{eqnarray}  However, since $\lim_{t \to N} r(t) = 1$, $\lim_{t \to N} \frac{u_l(t)}{u_l(N)} = 1$ for all $l \ge 0$, and $a_t(g)$ is given by some finite product of the form $r^{\sum_{s=1}^{n+1}|k(s)|} \prod_{s=1}^n \frac{u_{l(s)}(t)}{u_{l(s)}(N)}$, (\ref{pointwise_lim}) is obvious.  

To prove $(2)$, we need a better understanding of the quantities $a_t(g)$, ($g \in \F_2^+$).  Let $\chi_g$ denote the character of the irreducible corepresentation $U^g$.  From the fusion rules (\ref{eqn_fusion_unitary}) for $U_N^+$, it follows that the characters $\{\chi_g\}_{g \in \F_2^+}$ satisfy the relations \begin{eqnarray} \label{eqn_character_unitary}
\chi_g\chi_h = \sum_{\substack{\alpha, \beta, \sigma \in \F_2^+ \\ g = \alpha\sigma, \ h = \overline{\sigma}\beta}} \chi_{\alpha\beta} && (g,h \in \F_2^+). 
\end{eqnarray}
Furthermore,  since $U^{g_1} = z\boxtimes V^1$, we have
\begin{eqnarray} \label{IC_unitary}
\chi_{g_1} = z\chi_1 &\textrm{and}& \chi_{g_2} = \chi_{g_1}^* = \chi_1z^{-1},
\end{eqnarray} where $\{\chi_n\}_{n \ge 0}$ denote the irreducible characters of $O_N^+$, and $z = \id_\T$ is the generator of $C(\T)$ (which is $\ast$-free from $\{\chi_n\}_{n \ge 0}$).  Using (\ref{eqn_character_unitary}), the initial conditions (\ref{IC_unitary}), and the recursion (\ref{eqn_recursion_characters}) for the characters $\{\chi_n\}_{n \ge 0}$, it follows by induction on $|g|$ that for any $g \in \F_2^+$, there exist $k(1), \ldots, k(n) \in \N$ such that $|g| = k(1)+k(2)+ \ldots + k(n)$, $\epsilon(1), \epsilon(n+1) \in \{0, \pm 1\}$ and  $\epsilon(2), \ldots \epsilon(n) \in \{\pm 1\}$ such that $\chi_g$ is given by the free product
\begin{eqnarray} \label{eqn_char_again}
\chi_g = z^{\epsilon(1)}\chi_{k(1)}z^{\epsilon(2)}\chi_{k(2)} \ldots z^{\epsilon(n)}\chi_{k(n)}z^{\epsilon(n+1)}.
\end{eqnarray} 
Since $\Psi_t|_{L^2_g(U_N^+)} = a_t(g)\id_{L^2_g(U_N^+)}$ and $\chi_g \in L^2_g(U_N^+)$, we get from (\ref{eqn_char_again}) that \begin{eqnarray*}
a_t(g)\chi_g &=& \Psi_t\chi_g \\
&=& (P_{r(t)}*_{\textrm{red}}T_t)\big(z^{\epsilon(1)}\chi_{k(1)}z^{\epsilon(2)}\chi_{k(2)} \ldots z^{\epsilon(n)}\chi_{k(n)}z^{\epsilon(n+1)}\big) \\
&=& r(t)^{\sum_{s=1}^{n+1}|\epsilon(s)|} \prod_{s=1}^n \frac{u_{k(s)}(t)}{u_{k(s)}(N)}z^{\epsilon(1)}\chi_{k(1)}z^{\epsilon(2)}\chi_{k(2)} \ldots z^{\epsilon(n)}\chi_{k(n)}z^{\epsilon(n+1)} \\ 
&=& r(t)^{\sum_{s=1}^{n+1}|\epsilon(s)|} \prod_{s=1}^n \frac{u_{k(s)}(t)}{u_{k(s)}(N)} \chi_g \\
&=& r(t)^{|\epsilon(1)| + |\epsilon(2)| + n-1}\prod_{s=1}^n \frac{u_{k(s)}(t)}{u_{k(s)}(N)} \chi_g.
\end{eqnarray*}
Therefore $a_t(g)>0$, and using (\ref{eqn_q_cheby}) we have \begin{eqnarray*}
 a_t(g) &=&  r(t)^{|\epsilon(1)| + |\epsilon(2)| + n-1}\prod_{s=1}^n \frac{u_{k(s)}(t)}{u_{k(s)}(N)} \\
&\le&  r(t)^{n-1}\prod_{s=1}^n \frac{u_{k(s)}(t)}{u_{k(s)}(N)} \\ 
&=& r(t)^{n-1} \Big(\frac{q(t)}{q(N)}\Big)^{k(1) + \ldots + k(n)} \prod_{s=1}^n \frac{(1-q(t)^{-2k(s)-2})(1-q(N)^{-2})}{(1-q(N)^{-2k(s)-2})(1-q(t)^{-2})} \\
&=& r(t)^{-1} \Big(\frac{q(t)}{q(N)}\Big)^{k(1) + \ldots + k(n)} \prod_{s=1}^n \frac{1-q(t)^{-2k(s)-2}}{1-q(N)^{-2k(s)-2}} \\
&\le & \frac{1}{1-q(t)^{-2}} \Big(\frac{q(t)}{q(N)}\Big)^{k(1) + \ldots + k(n)} \\
&\le& \frac{1}{1-q(t_0)^{-2}} \Big(\frac{t}{N}\Big)^{k(1) + \ldots + k(n)} = C_{t_0} \Big(\frac{t}{N}\Big)^{|g|}.  
\end{eqnarray*}
\end{proof}

Using Proposition \ref{prop_decay_psi}, we can prove the MAP for $C(U_N^+)$ when $N \ge 3$.

\begin{thm} \label{thm_MAP_unitary}
For $N \ge 3$, $C(U_N^+)$ has the MAP.
\end{thm}

\begin{proof}
Since the proof is similar to the truncation argument for the orthogonal case, we only sketch it.  For each $m \in \N$, let \begin{eqnarray*} \Psi_{t,m} = \sum_{g \in \F_2^+: |g| \le m} a_t(g)p_g, & \textrm{and}& \tilde{\Psi}_{t,m} = \|\Psi_{t,m}\|_{B(C(U_N^+))}^{-1} \Psi_{t,m}.
\end{eqnarray*}
Then $\{\tilde{\Psi}_{t,m}\}_{t \in [t_0, N), m \in \N} \subset B(C(U_N^+))$ is net of finite rank contractions.  By Proposition \ref{prop_ultracontractive_unit} and Proposition \ref{prop_decay_psi}(2) 
\begin{eqnarray*}
\lim_{m \to \infty} \|\Psi_{t,m} - \Psi_t\|_{B(C(U_N^+))} \le \limsup_{m \to \infty} \sup_{n \ge m+1} (n+1)^2C_{t_0}\Big(\frac{t}{N}\Big)^n = 0,
\end{eqnarray*}
for every $t_0 \le t <N$.  Therefore $\lim_{m \to \infty} \|\Psi_{t,m}\|_{B(C(U_N^+))} = \|\Psi_t\|_{B(C(U_N^+))} = 1$, and consequently \begin{eqnarray} \label{norm_part}
\lim_{m \to \infty} \|\tilde{\Psi}_{t,m} - \Psi_t\|_{B(C(U_N^+))} = 0  && (t_0 \le t <N).
\end{eqnarray}  
Using (\ref{norm_part}) and Proposition \ref{prop_decay_psi}(1), it follows from arguments analogous to those in the proof of Theorem \ref{thm_MAP_orth}, that $\id_{C(U_N^+)}$ is contained in the point-norm closure of $\{\tilde{\Psi}_{t,m}\}_{t \in [t_0, N), m \in \N}$. 
\end{proof}

\subsubsection{The Case $N = 2$.}  

When $N = 2$, we can actually prove the following stronger approximation property for $C(U_2^+)$.

\begin{thm} \label{thm_CCAP_2}
$C(U_2^+)$ has the completely contractive approximation property (CCAP).  
\end{thm}

Recall that an operator space $Y$ has the CCAP if there is a net of finite rank complete contractions $\{\Phi_\lambda\}_{\lambda \in \Lambda} \subset CB(Y)$,  converging to $\id_Y$ in the point-norm topology.  A dual operator space $Y^*$ is said to have the \textit{weak$^\ast$ completely contractive approximation property} (weak$^\ast$-CCAP) if there is a net of $\sigma$-weakly continuous finite rank complete contractions  $\{\Phi_\lambda\}_{\lambda \in \Lambda} \subset CB(Y^*)$, converging to $\id_{Y^*}$ in the point $\sigma$-weak topology.

\begin{proof} 
In \cite[Theorem 5.14]{KrRu}, Kraus and Ruan showed that for a  compact quantum group $\G = (A, \Delta)$ of Kac type, $C(\G)$ has the CCAP if and only if $L^\infty(\G)$ has the weak$^\ast$-CCAP.  (In \cite{KrRu}, the authors state their results in the equivalent, dual framework of discrete Kac algebras.)  

However, as mentioned in the introduction, Banica has shown that $L^\infty(U_2^+) \cong L(\F_2)$ (see \cite{Ba}).  Since $L(\F_2)$ has the weak$^\ast$-CCAP and this property is an isomorphism invariant for von Neumann algebras (see \cite{CoHa}), $L^\infty(U_2^+)$ also has the weak$^\ast$-CCAP.  The theorem now follows from the previous paragraph.  
\end{proof}

\begin{rem} \label{rem_CCAP_quest} It would be interesting to know whether or not all of the C$^\ast$-algebras $C(\G)$ ($\G = O_N^+$ or $U_N^+$) have the CCAP, or perhaps just the completely bounded approximation property.  Since we have shown that these C$^\ast$-algebras always have the MAP, the answer seems likely to be yes.  (To our knowledge, there are no known examples of C$^\ast$-algebras with the MAP, but not the CCAP.)    We are currently investigating this problem.
\end{rem}

\begin{rem} \label{rem_bai_L1}
We conclude with the following Banach algebraic remark.  Let $\G$ be either $U_N^+$ for some $N \ge 2$, or $O_N^+$ for some $N \ge 3$.  Since $\G$ is not co-amenable, the Banach algebra $L^1(\G)$ fails to have a bounded approximate identity \cite[Theorem 3.1]{BeTu}.  However, our proofs of the MAP for $C(\G)$ can be used to show that $L^1(\G)$ has a central approximate identity, which is bounded in the \textit{(left) multiplier norm} on $L^1(\G)$: \begin{eqnarray*}
\|\omega\|_{M(L^1(\G))} := \sup\{\|\omega * \omega'\|_{L^1(\G)}: \|\omega'\|_{L^1(\G)}=1\} && (\omega \in L^1(\G)). 
\end{eqnarray*}  Indeed, when $\G = U_2^+$, $C(\G)$ has the CCAP by Theorem \ref{thm_CCAP_2}, and \cite[Theorems 5.14-5.15]{KrRu} show that $L^1(\G)$ has a central approximate identity which is actually bounded in the \textit{completely bounded} multiplier norm on $L^1(\G)$.  For all other cases, let $\{\Phi_t\}_{t \in 
\Lambda}$ be a net of finite rank contractions giving the MAP for $C(\G)$.  From the proofs of Theorems \ref{thm_MAP_orth} and \ref{thm_MAP_unitary}, we may assume that each $\Phi_t$ is a finite linear combination of the projections $\{p_\alpha\}_{\alpha \in \A}$ defined in Notation \ref{notat_L2_alpha}, say $\Phi_t = \sum_{\alpha \in \A} a_t(\alpha) p_\alpha$.  Identify $L^1(\G)$ with the closure of $L^\infty(\G)$ with respect to the norm $\|\omega\|_{L^1(\G)}:= \sup\{|h(\omega x)|: \|x\|_{L^\infty(\G)}=1\}$.  Then the net $\{\omega_t\}_{t \in \Lambda} \subset L^1(\G)$ defined by $\omega_t := \sum_{\alpha \in \A} a_t(\overline{\alpha}) d_\alpha \chi_\alpha$, is a central approximate identity for $L^1(\G)$ such that $\|\omega_t\|_{M(L^1(\G))} = \|\Phi_t\|_{B(C(\G))} = 1$ for all $t \in \Lambda$.  
The equality $\|\omega_t\|_{M(L^1(\G))} = \|\Phi_t\|_{B(C(\G))}$ follows from the simple verification that $\Phi_t = (m_{\omega_t})^*|_{C(\G)}$, where $m_{\omega_t}: \omega' \mapsto \omega_t*\omega'$ is the (left) multiplier induced by $\omega_t$.              
\end{rem}

\end{document}